\documentclass{article}[11pt]
\usepackage{graphicx} 
\usepackage[rgb]{xcolor}
\usepackage{tikz}
\usepackage{verbatim}
\usepackage{amsmath}
\usepackage{amsxtra}
\usepackage{amscd}
\usepackage{amsthm}
\usepackage{amsfonts}
\usepackage{amssymb}
\usepackage{eucal}
\usepackage{float}
\usepackage{mathrsfs}
\usepackage{graphicx}
\usepackage{stmaryrd}
\usepackage{tikz-cd}
\usetikzlibrary{arrows}
\usetikzlibrary{cd}
\usepackage{ytableau}

\usepackage{enumitem}
\setlist[enumerate]{label=(\arabic*), font=\normalfont}

\usepackage{latexsym,todonotes}

\usepackage[linktocpage=true]{hyperref}
\hypersetup{colorlinks,linkcolor=blue,urlcolor=orange,citecolor=blue}
\usepackage[all, cmtip]{xy}

\usepackage{stmaryrd}
\usepackage{color} 

\usepackage[alphabetic]{amsrefs}
\usepackage{todonotes}

\newtheorem{thm}{Theorem} [section]

\newtheorem{lem}[thm]{Lemma}
\newtheorem{prop}[thm]{Proposition}

\theoremstyle{definition}

\theoremstyle{remark}
\newtheorem{rem}[thm]{Remark}

\numberwithin{equation}{section}

\newcommand{\nc}{\newcommand}

\nc{\Ord}{\text{Ord}_v}

 \nc{\A}{\mathcal A} 
  \nc{\G}{\mathbb G} 
\nc{\Ainv}{\A^{\rm inv}}
\nc{\aA}{{}_\A}
\nc{\aAp}{{}_\A'}
\nc{\aff}{{}_\A\f}
\nc{\aL}{{}_\A L}
\nc{\aM}{{}_\A M}
\nc{\Bin}{B_i^{(n)}}
\nc{\dL}{{}^\omega L}
\nc{\Z}{{\mathbb Z}}
 \nc{\C}{{\mathbb C}}
 \nc{\N}{{\mathbb N}}
 \nc{\R}{{\mathbb R}}
  \nc{\I}{{\mathcal I}}
 \nc{\fZ}{{\mf Z}}
 \nc{\F}{{\mf F}}
 \nc{\Q}{\mathbb{Q}}
 \nc{\la}{\lambda}
 \nc{\ep}{\epsilon}
 \nc{\h}{\mathfrak h}
 \nc{\He}{\bold{H}}
 \nc{\htt}{\text{tr }}
 \nc{\n}{\mf n}
 \nc{\g}{{\mathfrak g}}
 \nc{\DG}{\widetilde{\mathfrak g}}
 \nc{\SG}{\breve{\mathfrak g}}
 \nc{\is}{{\mathbf i}}
 \nc{\V}{\mf V}
 \nc{\bi}{\bibitem}
 \nc{\E}{\mc E}
 \nc{\ba}{\tilde{\pa}}
 \nc{\half}{\frac{1}{2}}
 \nc{\hgt}{\text{ht}}
 \nc{\ka}{\kappa}
 \nc{\mc}{\mathcal}
 \nc{\mf}{\mathfrak} 
 \nc{\hf}{\frac{1}{2}}
\nc{\ov}{\overline}
\nc{\ul}{\underline}

\nc{\xx}{{\mf x}}
\nc{\id}{\text{id}}
\nc{\one}{\bold{1}}
\nc{\mfsl}{\mf{sl}}
\nc{\mfgl}{\mf{gl}}
\nc{\ti}[1]{\textit{#1}}
\nc{\Hom}{\text{Hom}}
\nc{\Cat}{\mathscr{C}}
\nc{\CatO}{\mathscr{O}}

\nc{\Tan}{\mathscr{T}}
\nc{\Umod}{\mathscr{U}}
\nc{\Func}{\mathscr{F}}
\nc{\Kh}{\text{Kh}}
\nc{\Khb}[1]{\llbracket #1 \rrbracket}

\nc{\ua}{\mf{u}}
\nc{\nb}{u}
\nc{\inv}{\theta}
\nc{\mA}{\mathcal{A}}

\nc{\Qq}{\Q(v)}
\nc{\U}{\bold{U}}
\nc{\uu}{\mathfrak{u}}
\nc{\Udot}{\dot{\U}}
\nc{\f}{\bold{f}}
\nc{\fprime}{\bold{'f}}
\nc{\B}{\bold{B}}
\nc{\Bdot}{\dot{\B}}
\nc{\Dupsilon}{\Upsilon^{\vartriangle}}

\newcommand{\frakd}{\mathfrak{d}}
\newcommand{\frakp}{\mathfrak{p}}
\newcommand{\calO}{\mathcal{O}}
\newcommand{\bbC}{\mathbb{C}}

\DeclareMathOperator{\tr}{tr}
\DeclareMathOperator{\Ker}{Ker}
\DeclareMathOperator{\PGL}{PGL}
\DeclareMathOperator{\SL}{SL}
\DeclareMathOperator{\Ind}{Ind}
\DeclareMathOperator{\Res}{Res}
\newcommand{\bbG}{\mathbb{G}}
\newcommand{\bbZ}{\mathbb{Z}}
\newcommand{\mapSJO}[1]{\Phi_{#1}}
\newcommand{\mapJOS}[1]{\Psi_{#1}}
\newcommand{\Nbent}[2]{N_{#1, #2}}
\DeclareMathOperator{\Gal}{Gal}
\DeclareMathOperator{\Lie}{Lie}
\newcommand{\sym}{\mathrm{sym}}
\newcommand{\asym}{\mathrm{asym}}
\newcommand{\pr}{\mathrm{pr}}

\nc{\ipsi}{\psi_{\imath}}
\nc{\Ui}{{\bold{U}^{\imath}}}
\nc{\uidot}{\dot{\mathfrak{u}}^{\imath}}
\nc{\Uidot}{\dot{\bold{U}}^{\imath}}
 \nc{\be}{e}
 \nc{\bff}{f}
 \nc{\bk}{k}
 \nc{\bt}{t}
 \nc{\bs}{\backslash}
 \nc{\BLambda}{{\Lambda_{\inv}}}
\nc{\Ktilde}{\widetilde{K}}
\nc{\bktilde}{\widetilde{k}}
\nc{\Yi}{Y^{w_0}}
\nc{\bunlambda}{\Lambda^\imath}

\nc{\ile}{\le_\imath}
\nc{\il}{<_{\imath}}

\nc{\etab}{\eta^{\bullet}}

\newcommand{\commentcustom}[1]{}

\begin{document}

\title{Double covers of tori for $p=2$ and genuine functions}
\author{L\'eo Gratien \and Nhat Hoang Le \and Kazuki Tokimoto}
\date{}
\maketitle

\begin{abstract}
In \cite{1907.05173}, Kaletha associates to a maximal torus \(S\) in a connected reductive group a canonical double cover \(S(F)_\pm\)
and introduces a genuine sign function \(a_S\) associated to a genuine character of \(S(F)_\pm\).
We give two constructions of this genuine sign function over a non-archimedean local field, including residual characteristic \(2\).
The first uses rootwise lifts to the double cover; the second, in characteristic zero, uses the exponential map.
Under explicit ramification conditions, both constructions yield genuine sign functions which, in odd residual characteristic, agree with the one constructed in \cite{1907.05173}.
We also illustrate the case of \(\mathrm{PGL}_2\).
A brief announcement of the first construction appears in \cite{1907.05173}*{Appendix D}.
\end{abstract}

\tableofcontents
\section{Introduction}

Let \(G\) be a connected reductive group over a non-archimedean local field \(F\)
of residual characteristic \(p\),
and let \(S\subset G\) be a maximal torus defined over \(F\).
In \cite{1907.05173}, Kaletha introduces 
a canonical topological double cover \(S(F)_{\pm}\) of \(S(F)\), which plays an important role in the harmonic analysis on \(G(F)\) and the local Langlands correspondence for \(G\).
On the Galois side, he associates an \(L\)-group to this double cover and constructs a canonical \(L\)-embedding into \({}^L G\), avoiding the choice of \(\chi\)-data required in the usual Langlands--Shelstad construction.
More directly relevant to the present paper, 
the double cover \(S(F)_{\pm}\) 
also provides a different point of view on
the character formula for regular (or more generally non-singular) supercuspidal representations of \(G(F)\)
and the explicit local Langlands correspondence for those representations
developed
in \cites{MR4013740, 1912.03274}.

The purpose of the present paper is to give two constructions
of a certain genuine function \(a_S\), which is an essential ingredient in the character formula,
without assuming that \(p\) is odd.
When \(p\neq 2\), a simpler construction is available in \cite{1907.05173}.
The constructions developed here are more involved, but have the advantage of not requiring the oddness assumption on \(p\).

Let us briefly recall 
the role of \(a_S\).
Extending the usual notion of genuine characters, we say that 
a function \(f \colon S(F)_{\pm} \to \bbC\) is \emph{genuine} if \(f(\epsilon \cdot x)=\epsilon \cdot f(x)\) for \(\epsilon \in \{\pm 1\}=\Ker (S(F)_{\pm}\twoheadrightarrow S(F))\) and \(x\in S(F)_{\pm}\).
The character formula associated to a genuine character \(\theta_{\pm}\) of \(S(F)_{\pm}\) in 
\cite{1907.05173}*{(4.3.2)}
is expressed 
in terms of a function on \(S(F)\),
and its value at \(s\in S(F)\) is of the following form:
\[
(\text{an explicit root of unity independent of \(\theta_{\pm}\)})
\cdot
\sum_{w\in N(S, G)(F)/S(F)} (a_S \cdot \theta_{\pm})(s^w),
\]
where \(N(S, G)\) is the normalizer of \(S\) in \(G\)
and \(a_S\) is a certain genuine function on \(S(F)_{\pm}\) constructed from \(\theta_{\pm}\).
By design, the product \(a_S\cdot \theta_{\pm}\) descends to \(S(F)\) and can be evaluated at \(s^w\in S(F)\).
Thus, the genuine function \(a_S\) is essential for the character formula.
It also appears in the conjectural stable character formula underlying the spectral characterization of
the local Langlands correspondence for supercuspidal \(L\)-parameters given in \cite{1907.05173}*{Conjecture 4.4.1, Remark 4.4.3}.
An announcement of the first construction appears in \cite{1907.05173}*{Appendix D}.
The present paper gives its complete details, develops a second construction, and studies both constructions in the case \(G=\PGL_2\).
\medskip

To discuss the difficulties in \(p=2\) and also indicate how we resolve them,
let us point out the two parts where the assumption \(p\neq 2\) matters in the simpler construction (see \S \ref{sec:gen_function} for more details).
It involves a set of mod-\(a\)-data in the sense of \cite{MR4013740}, which is a variant of \(a\)-data of Langlands--Shelstad \cite{MR909227}.
Recall that if \(\alpha\) is a symmetric root in the sense of \cite{MR909227},
then we have a tower \(F_{\alpha}/F_{\pm \alpha}/F\) of finite separable extensions with \(F_{\alpha}/F_{\pm \alpha}\) quadratic.
For each symmetric root \(\alpha\), one first obtains a mod-\(a\)-datum \(\overline{a}_{\alpha}\), and then chooses a trace-zero lift \(a_{\alpha}\in F_{\alpha}^0:=\Ker \tr_{F_{\alpha}/F_{\pm \alpha}}\) in order to define \(a_S\).
Two elementary facts are used here. First, such a lift can be obtained by antisymmetrizing an arbitrary lift and dividing by \(2\).
Second, the resulting genuine function is independent of the lift because the ratio of two lifts is a principal unit in \(F_{\pm \alpha}\),
and, since \(F_{\alpha}/F_{\pm \alpha}\) is tamely ramified when \(p\neq 2\), the quadratic character \(\kappa_{\alpha}\) attached to this extension is trivial on principal units.
Both arguments use essentially that \(p\neq2\).
This part is reviewed in \S \ref{sec:modatoas}.

To treat the case \(p=2\), we 
use the additive auto-duality of \(F_{\alpha}\) to reformulate the existence of a trace-zero lift as a condition on the restriction of a character
(Proposition \ref{prop:liftofa}),
which is verified for the characters arising in both constructions,
and we ``refine'' the datum \(\overline{a}_{\alpha}\) to ensure that the ratio of two lifts lies in the kernel of \(\kappa_{\alpha}\) (Proposition \ref{prop:wildmoda}).

The other part where the oddness assumption on \(p\) is used is the step to obtain a set of mod-\(a\)-data from the genuine character \(\theta_{\pm}\).
Fix a non-trivial character \(\Lambda\) of \(F\).
One takes the square \(\theta_{\pm}^2\) to descend it to \(S(F)\) and for each symmetric root \(\alpha\) considers the composite \(\theta_{\pm}^2 \circ N_{F_{\alpha}/F} \circ \alpha^{\vee}\), which is a character of \(F_{\alpha}^{\times}\).
If the depth \(r_{\theta_{\pm}, \alpha}\) of this character is positive, then one defines \(\overline{a}_{\alpha}\) by requiring the following equality
\[
\theta_{\pm}^2(N_{F_{\alpha}/F}(\alpha^{\vee}(z+1)))=\Lambda (\tr_{F_{\alpha}/F}(2\overline{a}_{\alpha}z))
\]
to hold for any \(z\in [F_{\alpha}]_{r_{\theta_{\pm}, \alpha}}\),
where \([F_{\alpha}]_{r}:=\{z\in F_{\alpha} \mid v(z)\geq r\}\) and \(v\) is the extension to \(F_{\alpha}\) of the normalized valuation of \(F\).
To summarize,
one first extracts a character of \(F_{\alpha}^{\times}\) from \(\theta_{\pm}\)
and then uses the Moy--Prasad isomorphism and the additive duality of \(F_{\alpha}\).
Note that the factor \(2\) on the right-hand side 
corresponds to the square on the left-hand side;
this squaring-and-doubling procedure is problematic when \(p=2\).
This part is reviewed in \S \ref{sec:gctoas}.

\medskip 

The key idea in both constructions is therefore to extract ``the \(\alpha\)-part'' of \(\theta_{\pm}\) directly, without squaring \(\theta_{\pm}\).
The two methods differ in how this is done.
The first method, sketched in \cite{1907.05173}*{Appendix D}, involves constructing a homomorphism \(F_{\alpha}^{\times} \to S(F)_{\pm}\) that lifts \(N_{F_{\alpha}/F} \circ \alpha^{\vee} \colon F_{\alpha}^{\times} \to S(F)\) and factors through \(F_{\alpha}^{\times} \twoheadrightarrow F_{\alpha}^{\times}/N_{F_{\alpha}/F_{\pm \alpha}}(F_{\alpha}^{\times})\).
Pulling back \(\theta_{\pm}\) via this homomorphism yields a character of \(F_{\alpha}^{\times}\), with which we can define \(\overline{a}_{\alpha}\) in a similar way.
The second method assumes that \(F\) is of characteristic zero and uses the exponential map of \(S(F)_{\pm}\) studied in \cite{1907.05173}*{\S 3.7} together with the differential maps of \(N_{F_{\alpha}/F}\) and \(\alpha^{\vee}\) to pull back \(\theta_{\pm}\) and obtain a character of a suitable open subgroup of \(F_{\alpha}\) where the relevant exponential maps converge.
Then we can define \(\overline{a}_{\alpha}\) via the additive duality of \(F_{\alpha}\).

In fact, due to the ``refinement'' of the mod-\(a\)-data, 
both methods require certain assumptions, denoted by \({\bf (Ram)}\) and \({\bf (Exp)}\), obtained by imposing the rootwise conditions \({\bf (Ram_{\alpha})}\) and \({\bf (Exp_{\alpha})}\), respectively, for every symmetric root \(\alpha\).
Roughly speaking, these conditions say that ``the \(\alpha\)-part'' of \(\theta_\pm\) is sufficiently more wildly ramified than the extension \(F_\alpha/F_{\pm\alpha}\).
We feel that assumptions of this sort are natural.
When \(p\) is odd, {\bf (Ram\(_\alpha\))} is equivalent to the positive-depth condition in \cite{1907.05173}, which we write as {\bf (PD\(_\alpha\))} in \S \ref{sec:gctoas}, and hence {\bf (Ram)} is equivalent to {\bf (PD)}.
When \(p=2\), {\bf (Ram)} agrees with condition \((*)\) in \cite{1907.05173}*{Theorem D.3.1}.
In odd residual characteristic, by contrast, \((*)\) is vacuous whereas {\bf (Ram)} imposes the positive-depth condition; see Remark \ref{rem:Ramandast}.
The rootwise condition {\bf (Ram\(_\alpha\))} also has another interpretation. As discussed in Remark \ref{rem:epsconstr}, it allows us to apply Deligne's twist formula \cite{MR349635} to the relevant characters of \(F_\alpha^\times\). When \(F_\alpha/F_{\pm\alpha}\) is ramified, the inequality in {\bf (Ram\(_\alpha\))} is precisely the ramification condition required for the twist formula, which expresses the \(\alpha\)-factor of \(a_S\) in terms of local epsilon factors.

\medskip

We also illustrate our two constructions by discussing the simplest example where \(G=\PGL_2\) and \(S\subset G\) is an elliptic maximal torus.
If \(E/F\) is the separable quadratic extension splitting \(S\), then a genuine character \(\theta_{\pm}\) can be identified with a character of \(E^{\times}\) trivial on \(N_{E/F}(E^{\times})\) but non-trivial on \(F^{\times}\).
We make explicit the assumptions {\bf (Ram\(_{\alpha}\))} and {\bf (Exp\(_{\alpha}\))} in this case, and in particular see that neither condition implies the other in general if \(E/F\) is wildly ramified.

\medskip 

This paper is organized as follows.
In \S \ref{sec:background}, we recall the relevant part of the theory in \cite{1907.05173}.
In \S \ref{sec:dc} we review the construction of the double cover \(S(F)_{\pm}\) of \(S(F)\) as well as a way to express general elements of \(S(F)_{\pm}\),
which is then used in \S \ref{sec:gen_function} to describe the construction of the genuine function \(a_S\) when \(p\neq 2\).
In \S \ref{sec:refinedmoda}, we give a few propositions used in our constructions of \(a_S\).
The first two propositions are used in both constructions to handle the first difficulty mentioned above, and the last proposition is used in variant constructions of the first method given in \S \ref{sec:relspl}.
In \S \ref{sec:constr}, we explain our two constructions of \(a_S\).
In \S \ref{sec:liftJOS}, we construct a certain homomorphism which essentially yields the homomorphism \(F_{\alpha}^{\times} \to S(F)_{\pm}\) mentioned above, and in \S \ref{sec:method1} we give the first construction of \(a_S\) using this homomorphism.
In \S \ref{sec:method2}, we first recall the exponential map of \(S(F)_{\pm}\) and then give the second construction of \(a_S\) using it.
In \S \ref{sec:PGL2}, we illustrate our two constructions in the case \(G=\PGL_2\).
Roughly speaking, our announcement in \cite{1907.05173}*{Appendix D} is a compressed summary of \S \ref{sec:refinedmoda}, \S \ref{sec:liftJOS}, and \S \ref{sec:method1} of the present paper, and the present paper also contains the second construction in \S \ref{sec:method2} and the example in \S \ref{sec:PGL2}.
For the reader's convenience, we indicate the corresponding propositions and theorem in \cite{1907.05173}*{Appendix D} at the relevant places in the present paper.

\subsection*{Acknowledgments}
This work grew out of a project at the Arizona Winter School 2025.
We thank Tasho Kaletha for his lectures and for introducing the questions considered here, Jialiang Zou for his help during the project, and the organizers of the Arizona Winter School for providing a stimulating setting in which this work began.
We also thank Tasho Kaletha for giving us the opportunity to include an announcement of this work as Appendix D of \cite{1907.05173}.
K.T. was supported by JSPS KAKENHI Grant Number JP19K14503.

\subsection*{Notation}
Throughout the paper, \(F\) denotes a non-archimedean local field of residue characteristic \(p\), and \(F^{\mathrm{s}}\) is a fixed separable closure. For every finite extension \(K/F\) contained in \(F^{\mathrm{s}}\), we write \(\calO_K\) for its ring of integers, \(\frakp_K\) for the maximal ideal, and \(k_K=\calO_K/\frakp_K\) for the residue field. The normalized additive valuation is denoted by \(v_K\), so that \(v_K(K^{\times})=\bbZ\). For \(n\in\bbZ\), we write \(\frakp_K^n=\{x\in K\mid v_K(x)\geq n\}\), while
\[
U_K=\calO_K^{\times},\qquad U_K^n=1+\frakp_K^n\quad(n\geq1).
\]
The \emph{level} of a character of \(K^{\times}\) is the smallest integer \(l\geq 0\) for which it is trivial on \(U_{K}^{l+1}\).

If \(L/K\) is a finite separable extension, \(\tr_{L/K}\) and \(N_{L/K}\) denote the field trace and norm. When \(L/K\) is Galois, \(\Gamma_{L/K}=\Gal(L/K)\). We write \(\Gamma_K=\Gal(F^{\mathrm{s}}/K)\) and \(W_K\) for the absolute Galois and Weil groups of \(K\). The symbol \(\Res_{L/K}\) denotes restriction of scalars.

For an algebraic group \(H\) over \(K\), let \(\Lie(H)\) be its Lie algebra. If \(H\) is commutative and \(L/K\) is finite, we use \(N_{L/K}:H(L)\to H(K)\) for the group-theoretic norm and \(\tr_{L/K}:\Lie(H)(L)\to\Lie(H)(K)\) for its differential; the intended meaning is clear from the source and target. For a torus \(S\), let \(X^{\ast}(S)\) and \(X_{\ast}(S)\) denote its absolute character and cocharacter groups.

\section{Background on double covers and the genuine function} \label{sec:background}
Let \(G\) be a connected reductive group over \(F\)
and \(S\) a maximal torus in \(G\).
Here we briefly recall the relevant part of the theory developed in \cite{1907.05173}: the constructions of the double cover \(S(F)_{\pm}\) of \(S(F)\) and the genuine sign function \(a_S\) for odd \(p\).
For more complete details, we refer the reader to \cite{1907.05173}.

\subsection{Double covers of tori} \label{sec:dc}
\subsubsection*{Root data and associated fields}
Let \(X^{\ast}(S)\) be the absolute character group of \(S\) and \(R\subset X^{\ast}(S)\) be the set of absolute roots of \(G\) with respect to \(S\).
Thus, \(R\subset X^{\ast}(S)\) is stable under the natural action of \(\Sigma:=\Gamma_F \times \{\pm 1\}\).
In \cite{1907.05173}, one works more generally with a finite set \(R\) with a suitable action of \(\Sigma\) and a \(\Sigma\)-equivariant map \(R\to X^{\ast}(S)\).
However, for our purpose, we only need the case where \(R\) is the set of absolute roots of \(G\) with respect to \(S\).

Following \cite{MR909227}, a \(\Sigma\)-orbit \(O\subset R\) is said to be \emph{symmetric} if \(O\) is a single \(\Gamma_F\)-orbit and \emph{asymmetric} otherwise.
Similarly, a root \(\alpha \in R\) is said to be \emph{symmetric} or \emph{asymmetric} according to its \(\Sigma\)-orbit.
For \(\alpha \in R\), let \(\Gamma_{\alpha}\) and \(\Gamma_{\pm \alpha}\) be the stabilizers of \(\alpha\) and \(\{\pm \alpha\}\), and let \(F_{\alpha}\) and \(F_{\pm \alpha}\) be their respective fixed fields.
Then \(\alpha \in R\) is symmetric (resp.\ asymmetric) if and only if \([F_{\alpha}:F_{\pm \alpha}]=2\) (resp.\ \(F_{\alpha}=F_{\pm \alpha}\)).
The corresponding subsets of \(R\) are denoted by \(R_{\sym}\) and \(R_{\asym}\).

For \(\alpha \in R_{\sym}\), we have a quadratic separable extension \(F_{\alpha}/F_{\pm \alpha}\).
We write \(\tau_{\alpha}\) for the generator of the Galois group \(\Gamma_{F_{\alpha}/F_{\pm \alpha}}:=\Gal(F_{\alpha}/F_{\pm \alpha})\)
and \(\kappa_{\alpha}\) for the isomorphism \(F_{\pm \alpha}^{\times}/N_{F_{\alpha}/F_{\pm \alpha}}(F_{\alpha}^{\times}) \stackrel{\sim}{\to} \{\pm 1\}\) provided by the local class field theory.
Moreover, we put 
\[
F_{\alpha}^1:=\Ker N_{F_{\alpha}/F_{\pm \alpha}}, \quad
F_{\alpha}^0:=\Ker \tr_{F_{\alpha}/F_{\pm \alpha}}.
\]
For later convenience, we introduce the notation \(\widetilde{N}_{\alpha}\) for the natural surjection 
\[
\widetilde{N}_{\alpha} \colon F_{\alpha}^{\times} \to F_{\alpha}^1, \quad x \mapsto x/\tau_{\alpha}(x),
\]
whose kernel is \(F_{\pm \alpha}^{\times}\).

\subsubsection*{The auxiliary tori \(J_O\)}
We recall from \cite[\S3.1]{1907.05173} the construction of the double cover \(S(F)_{\pm}\) associated to \(R\subset X^{\ast}(S)\).
Let \(O\subset R\) be a \(\Sigma\)-orbit and choose \(\alpha\in O\). Define a torus \(J_{\alpha}\) by
\[
\begin{aligned}
J_{\alpha}
&:=
\begin{cases}
\Res_{F_{\alpha}/F}\G_m
& \text{if \(\alpha \in R_{\asym}\)}\\[2pt]
\Res_{F_{\pm\alpha}/F}
\Res^1_{F_{\alpha}/F_{\pm\alpha}}\G_m
& \text{if \(\alpha \in R_{\sym}\)},
\end{cases}
\end{aligned}
\]
where 
\(\Res_{F_{\alpha}/F_{\pm \alpha}}^1\G_m\) is the unique one-dimensional anisotropic torus defined over \(F_{\pm \alpha}\) and split over \(F_{\alpha}\).
Thus, \(J_{\alpha}(F)\) equals \(F_{\alpha}^{\times}\) or \(F_{\alpha}^1\) according to whether \(\alpha\) is asymmetric or symmetric.

In both cases, there is a natural morphism \(S\to J_{\alpha}\) over \(F\). For \(\alpha,\beta\in O\), there is a natural \(F\)-isomorphism \(J_{\alpha}\simeq J_{\beta}\), compatible with these morphisms. We define \(J_O\) as the inverse limit of the \(J_{\alpha}\) for \(\alpha\in O\), and it is equipped with a natural morphism \(S\to J_O\).

\subsubsection*{The orbitwise covers \(J_O(F)_{\pm}\)}
We next define a double cover \(J_O(F)_{\pm}\to J_O(F)\). If \(O\) is asymmetric, it is the split cover \(J_O(F)_{\pm}=J_O(F)\times\{\pm1\}\).

Suppose that \(O\) is symmetric and choose \(\alpha\in O\).
Then we define \(J_{\alpha}(F)_{\pm}:=F_{\alpha}^{\times}/N_{F_{\alpha}/F_{\pm \alpha}}(F_{\alpha}^{\times})\).
As \(J_{\alpha}(F)=F_{\alpha}^1\), 
the homomorphism \(\widetilde{N}_{\alpha}\) induces the following surjection, for which we use the same notation:
\[
\widetilde{N}_{\alpha} \colon
J_{\alpha}(F)_{\pm}=F_{\alpha}^{\times}/N_{F_{\alpha}/F_{\pm \alpha}}(F_{\alpha}^{\times}) \twoheadrightarrow J_{\alpha}(F)=F_{\alpha}^1, 
\]
whose kernel is \(F_{\pm \alpha}/N_{F_{\alpha}/F_{\pm \alpha}}(F_{\alpha}^{\times}) \stackrel{\kappa_{\alpha}}{\simeq} \{\pm 1\}\).

For \(\alpha,\beta\in O\), the transition isomorphism \(J_{\alpha}\simeq J_{\beta}\) lifts naturally to an isomorphism of the covers.
Thus, we define \(J_O(F)_{\pm}\) as the inverse limit of \(J_{\alpha}(F)_{\pm}\) over all \(\alpha \in O\) and it is equipped with a natural surjection \(J_O(F)_{\pm} \twoheadrightarrow J_O(F)\) whose kernel is isomorphic to \(\{\pm 1\}\).

\subsubsection*{The double cover of \(S(F)\)}
We define the \emph{big cover} \(S(F)_{\pm \pm}\) of \(S(F)\) as the fibered product
\[
S(F)_{\pm\pm}
:=
S(F)\times_{\prod_O J_O(F)}\prod_O J_O(F)_{\pm},
\]
where the products are taken over all \(\Sigma\)-orbits \(O\subset R\).
It is an extension of \(S(F)\) by \(\prod_{O\subset R} \{\pm 1\}\).
Then the double cover \(S(F)_{\pm}\) of \(S(F)\) is defined as the pushout of \(S(F)_{\pm \pm}\) along the multiplication map \(\prod_{O\subset R} \{\pm 1\} \to \{\pm 1\}\).

\subsubsection*{Expressions for elements of the cover}
We discuss a way to express a general element of \(S(F)_{\pm}\).
We use a slightly more redundant description than \cite[Remark 3.1.3]{1907.05173}.

For an asymmetric \(\Sigma\)-orbit \(O\subset R_{\asym}\), an element of \(J_O(F)_{\pm}\) is expressed as \((\gamma_{O}, \epsilon_O)\) with \(\gamma_{O} \in J_O(F)\) and \(\epsilon_O \in \{\pm 1\}\).

For a symmetric \(\Sigma\)-orbit \(O\subset R_{\sym}\), an element of \(J_O(F)\) is expressed as a tuple \((\gamma_{\alpha})_{\alpha \in O}\) such that \(\gamma_{\alpha} \in J_{\alpha}(F)=F_{\alpha}^1\) and \(\gamma_{\sigma \alpha} = \sigma(\gamma_{\alpha})\) for any \(\sigma \in \Gamma_F\),
while an element of \(J_O(F)_{\pm}\) is expressed as a tuple \((\overline{\delta}_{\alpha})_{\alpha \in O}\) such that \(\overline{\delta}_{\alpha} \in J_{\alpha}(F)_{\pm}=F_{\alpha}^{\times}/N_{F_{\alpha}/F_{\pm \alpha}}(F_{\alpha}^{\times})\) and \(\overline{\delta}_{\sigma \alpha}=\sigma(\overline{\delta}_{\alpha})\) for any \(\sigma \in \Gamma_F\).
Following \cite{1907.05173}, we usually choose a representative \(\delta_{\alpha} \in F_{\alpha}^{\times}\) of \(\overline{\delta}_{\alpha} \in J_{\alpha}(F)_{\pm}\) for each \(\alpha \in O\) such that \(\sigma(\delta_{\alpha}) = \delta_{\sigma \alpha}\) for \(\sigma \in \Gamma_F\) and simply use \((\delta_{\alpha})_{\alpha \in O}\) to represent an element of \(J_O(F)_{\pm}\).

Thus, in view of the condition given by the fibered product, an element of \(S(F)_{\pm \pm}\) is naturally expressed as a tuple 
\begin{equation}\label{eq:eltofdc}
(\gamma, (\delta_{\alpha})_{\alpha}, (\epsilon_O)_{O})=
(\gamma, (\delta_{\alpha})_{\alpha \in R_{\sym}}, (\epsilon_O)_{O\subset R_{\asym}})
\end{equation}
with \(\gamma \in S(F)\), \(\delta_{\alpha} \in F_{\alpha}^{\times}\), and \(\epsilon_O \in \{\pm 1\}\)
satisfying
\[
\delta_{\sigma \alpha}=\sigma(\delta_{\alpha}),
\qquad 
\alpha(\gamma)=\frac{\delta_{\alpha}}{\delta_{-\alpha}}
\]
for \(\sigma \in \Gamma_F\) and \(\alpha \in R_{\sym}\).
Since \(S(F)_{\pm}\) is a quotient of \(S(F)_{\pm \pm}\), an element of \(S(F)_{\pm}\) can in turn be represented by an element of \(S(F)_{\pm \pm}\).
We use this expression to define a genuine sign function \(a_S\) on \(S(F)_{\pm}\) in \S \ref{sec:gen_function}.

\subsubsection*{Genuine characters and \(\chi\)-data}
As usual, a character of \(S(F)_{\pm}\) is said to be \emph{genuine} if it does not factor through the projection \(S(F)_{\pm} \twoheadrightarrow S(F)\).
Let us also recall the notion of \(\chi\)-data introduced in \cite{MR909227}*{\S2.5} and briefly indicate its close connection with the theory of double covers.
A set of \(\chi\)-data \(\chi=\left(\chi_{\alpha}\right)_{\alpha \in R}\) consists of a character \(\chi_{\alpha}\colon F_{\alpha}^{\times} \to \C^{\times}\) for each \(\alpha \in R\) such that \(\chi_{\sigma \alpha}=\chi_{\alpha} \circ \sigma^{-1}\) and \(\chi_{-\alpha}=\chi_{\alpha}^{-1}\) for any \(\sigma \in \Gamma_F\) and any \(\alpha \in R\),
and \(\chi_{\alpha}|_{F_{\pm \alpha}^{\times}}=\kappa_{\alpha}\) for any \(\alpha \in R_{\sym}\).
The last condition implies that \(\chi_{\alpha}\) induces a genuine character of \(J_{\alpha}(F)_{\pm}\) for any \(\alpha \in R_{\sym}\).
This observation is developed further in \cite{1907.05173}*{\S3.2}, where it is shown that a set of \(\chi\)-data \(\chi\) gives rise to a genuine character \(\chi_S\) of \(S(F)_{\pm}\).

If \(\alpha\) is symmetric and \(\chi_{\alpha}|_U\) is valued in \(\{\pm 1\} \subset \bbC^{\times}\) for some open subgroup \(U\subset F_{\alpha}^{\times}\), then it induces a splitting of the cover over \(\widetilde{N}_{\alpha}(U) \subset J_{\alpha}(F)\).
More precisely, putting
\(J_{\alpha, U}(F)_{\pm}:=\widetilde{N}_{\alpha}^{-1}(\widetilde{N}_{\alpha}(U))=F_{\pm \alpha}^{\times}U/N_{F_{\alpha}/F_{\pm \alpha}}(F_{\alpha}^{\times})\), we have an isomorphism
\[
J_{\alpha, U}(F)_{\pm} \xrightarrow{\sim} \widetilde{N}_{\alpha}(U) \times \{\pm 1\}, \quad \gamma \mapsto (\widetilde{N}_{\alpha}(\gamma), \chi_{\alpha}(\gamma)).
\]
If \(\theta_{\pm}\) is a genuine character of \(J_{\alpha}(F)_{\pm}\),
then the restriction to \(J_{\alpha, U}(F)_{\pm}\) descends to a character of \(\widetilde{N}_{\alpha}(U)\) via the above splitting.
One can check that it is nothing but the restriction of the character of \(J_{\alpha}(F)\) obtained by descending the non-genuine character \(\theta_{\pm} \cdot \chi_{\alpha}^{-1}\).
Note that \(\chi_{\alpha}^{-1}=\chi_{\alpha}\) on \(J_{\alpha, U}(F)_{\pm}\). 
We will use this fact in \S \ref{sec:relspl}.
See \cite{1907.05173}*{\S 3.6} for the general relation between \(\chi\)-data and splittings of the double cover.

\subsection{Genuine function associated to genuine characters for odd \(p\)} \label{sec:gen_function}

Let \(\theta_{\pm}\) be a genuine character of \(S(F)_{\pm}\).
In \cite{1907.05173}*{\S4.3}, a possible formula is given for the restriction to (a certain subset of) \(S(F)\subset G(F)\) of the Harish-Chandra character of a supercuspidal representation of \(G(F)\) expected to be associated to \((S, \theta_{\pm})\).
Under the relevant additional hypotheses of \cite{MR4013740}, such representations are known to exist when \(p\) is odd.
One of the ingredients of the formula is a certain sign function \(a_S\colon S(F)_{\pm} \to \{\pm 1\}\).
It is \emph{genuine} in the sense of \cite{1907.05173}*{Fact 3.5.1}, i.e.\ \(a_S(\epsilon \cdot \tilde{\gamma}) = \epsilon \cdot a_S(\tilde{\gamma})\) for any \(\tilde{\gamma} \in S(F)_{\pm}\) and \(\epsilon \in \{\pm 1\}=\Ker(S(F)_{\pm} \twoheadrightarrow S(F))\).

The character formula in \cite{1907.05173}*{\S4.3} is stated without assumptions on the base field \(F\).
For \(p=2\), however, the construction of \(a_S\) needed for the formula is announced in \cite{1907.05173}*{Appendix D} and developed in the present paper.
Here we recall the construction for odd \(p\) in three parts,
highlighting how the ingredients are related and where the
assumption that \(p\) is odd enters.

\subsubsection{From \(a\)-data to a genuine function} \label{sec:atoas}
We first recall the construction of a genuine function from \(a\)-data
given in \cite{1907.05173}*{\S3.5}. This construction does not require
any assumption on \(p\).
Recall from \cite{MR909227}*{\S2.2} that a set of \(a\)-data \(\left(a_{\alpha}\right)_{\alpha \in R}\) consists of an element \(a_{\alpha} \in F_{\alpha}^{\times}\) for each \(\alpha \in R\) such that \(a_{\sigma \alpha}=\sigma(a_{\alpha})\) and \(a_{-\alpha}=-a_{\alpha}\) for any \(\sigma \in \Gamma_F\) and \(\alpha \in R\). 
Given a set of \(a\)-data, we first define a sign function on the big cover \(S(F)_{\pm \pm}\) as follows (by slight abuse of notation we write it also as \(a_S\)):
\begin{align} \label{eq:aSfroma}
a_S \colon S(F)_{\pm \pm} &\to \{\pm 1\} \\
(\gamma, (\delta_{\alpha})_{\alpha}, (\epsilon_O)_{O})
&\mapsto \prod_{\substack{\alpha \in R_{\sym}/\Gamma \\ \alpha(\gamma)\neq 1}} \kappa_{\alpha}\left(\frac{\delta_{\alpha}-\delta_{-\alpha}}{a_{\alpha}}\right)
\cdot \prod_{\substack{\alpha \in R_{\sym}/\Gamma \\ \alpha(\gamma)=1}} \kappa_{\alpha}(\delta_{\alpha})
\cdot \prod_{O\subset R_{\asym}} \epsilon_O. \notag
\end{align}
Here, an element of \(S(F)_{\pm \pm}\) is expressed as in \eqref{eq:eltofdc}.
In the first product, as the elements \(\delta_{\alpha}-\delta_{-\alpha}=\delta_{\alpha}-\tau_{\alpha}(\delta_{\alpha})\) and \(a_{\alpha}\) both lie in \(F_{\alpha}^0\setminus \{0\}\), the quotient \(\frac{\delta_{\alpha}-\delta_{-\alpha}}{a_{\alpha}}\) lies in \(F_{\pm \alpha}^{\times}\) and its image under \(\kappa_{\alpha}\) only depends on the image of \(\delta_{\alpha}\in F_{\alpha}^{\times}\) in \(J_{\alpha}(F)_{\pm}\).
Similarly, in the second product, the element \(\delta_{\alpha}\) lies in \(F_{\pm \alpha}^{\times}\) and its image under \(\kappa_{\alpha}\) only depends on the image of \(\delta_{\alpha}\) in \(J_{\alpha}(F)_{\pm}\).

Moreover, it is easy to check that the above function descends to a genuine function \(a_S\colon S(F)_{\pm} \to \{\pm 1\}\); in other words,
we have 
\[
a_S((\epsilon_O)_O \cdot \tilde{\gamma}) = (\prod_{O\subset R} \epsilon_O) \cdot a_S(\tilde{\gamma})
\]
for any \(\tilde{\gamma} \in S(F)_{\pm \pm}\) and \((\epsilon_O)_O \in \prod_{O\subset R}\{\pm 1\}=\Ker(S(F)_{\pm \pm} \twoheadrightarrow S(F))\).

Note that the formula \eqref{eq:aSfroma} depends on \(a_{\alpha}\) only for \(\alpha \in R_{\sym}\). Thus, by slight abuse of terminology, we may specify only the \(a\)-data for symmetric roots when constructing the associated genuine function in what follows.

\subsubsection{From mod-\(a\)-data to a genuine function} \label{sec:modatoas}
We next recall the variant of the preceding construction discussed in
\cite{1907.05173}*{Remark 3.5.2}, in which one starts with
mod-\(a\)-data rather than \(a\)-data.
The assumption that \(p\) is odd enters at two points in passing from mod-\(a\)-data to a genuine function.

Recall from \cite{MR4013740}*{Definition 4.6.8} that a set of mod-\(a\)-data consists of a real number \(r_{\alpha}\) and a non-zero element \(\overline{a}_{\alpha} \in [F_{\alpha}]_{r_{\alpha}}/[F_{\alpha}]_{r_{\alpha}+}\) for each \(\alpha \in R\) such that \(r_{\sigma \alpha}=r_{\alpha}=r_{-\alpha}\), \(\overline{a}_{\sigma \alpha}=\sigma(\overline{a}_{\alpha})\) and \(\overline{a}_{-\alpha}=-\overline{a}_{\alpha}\) for any \(\sigma \in \Gamma_F\) and \(\alpha \in R\).
Here, we put \([F_{\alpha}]_{r_{\alpha}}:=\{x \in F_{\alpha} \mid v(x) \geq r_{\alpha}\}\) and \([F_{\alpha}]_{r_{\alpha}+}:=\{x \in F_{\alpha} \mid v(x) > r_{\alpha}\}\), where \(v\) is the normalized valuation on \(F\) extended to \(F_{\alpha}\).

As discussed in \cite{1907.05173}*{Remark 3.5.2}, if \(p\) is odd, mod-\(a\)-data can be lifted to \(a\)-data and the genuine function obtained from \S \ref{sec:atoas} is independent of the choice of lift.
Let us briefly recall the argument.
For the existence of \(a\)-data lifting the given mod-\(a\)-data, it suffices to show that \(\overline{a}_{\alpha}\) can be lifted to an element \(a_{\alpha} \in F_{\alpha}^0\) for each \(\alpha \in R_{\sym}\).
This can be achieved by taking an arbitrary lift \(\dot{a}_{\alpha} \in [F_{\alpha}]_{r_{\alpha}}\) and setting \(a_{\alpha}=\frac{1}{2}(\dot{a}_{\alpha} - \tau_{\alpha}(\dot{a}_{\alpha})) \in F_{\alpha}^0\) because \(p\neq 2\).

For two lifts \(\left(a_{\alpha}\right)_{\alpha}\) and \(\left(a'_{\alpha}\right)_{\alpha}\) of the same mod-\(a\)-data, we have \(a'_{\alpha}a_{\alpha}^{-1} \in U_{F_{\alpha}}^1\cap F_{\pm \alpha}=U_{F_{\pm \alpha}}^1\) for any \(\alpha \in R_{\sym}\), but \(\kappa_{\alpha}\) is trivial on \(U_{F_{\pm \alpha}}^1\) since \(p\neq 2\) and thus \(F_{\alpha}/F_{\pm \alpha}\) is tamely ramified.
This shows that the two lifts define the same genuine function.

\subsubsection{From a genuine character to a genuine function} \label{sec:gctoas}
Finally, we recall from \cite{1907.05173}*{\S4.3} how a genuine
character \(\theta_\pm\) gives rise to a genuine function when
\(p\) is odd.
At positive depth, the construction first extracts mod-\(a\)-data
from \(\theta_{\pm}\) and then applies \S \ref{sec:modatoas}; the assumption that \(p\) is odd enters in both steps.
The depth-zero case can be treated directly.

Fix a non-trivial character \(\Lambda\) of \(F\).
The square \(\theta_{\pm}^2\) of \(\theta_{\pm}\) descends to a character of \(S(F)\) and we can consider the character \(\theta_{\pm}^2 \circ N_{F_{\alpha}/F} \circ \alpha^{\vee}\) of \(F_{\alpha}^{\times}\).
Let \(r_{\theta_{\pm}, \alpha}\) be the depth of this character and let \(r_{\Lambda, \alpha}\) be the depth of the character \(\Lambda \circ \tr_{F_{\alpha}/F}\).
For \(\alpha \in R_{\sym}\) with \(r_{\theta_{\pm}, \alpha}>0\), define \(r_{\alpha}:=r_{\Lambda, \alpha}-r_{\theta_{\pm}, \alpha}\).
For each \(\alpha \in R_{\sym}\),
we assume the following positive-depth condition:
\begin{description}
	\item[(PD\(_{\alpha}\))] if \(F_{\alpha}/F_{\pm \alpha}\) is ramified, then \(r_{\theta_{\pm}, \alpha}>0\).
\end{description}
We write {\bf (PD)} for the condition that {\bf (PD\(_{\alpha}\))} holds for every \(\alpha \in R_{\sym}\).

For \(\alpha \in R_{\sym}\) such that \(r_{\theta_{\pm}, \alpha} >0\), there exists a unique non-zero element \(\overline{a}_{\alpha} \in [F_{\alpha}]_{r_{\alpha}}/[F_{\alpha}]_{r_{\alpha}+}\) such that
\[
\theta_{\pm}^2(N_{F_{\alpha}/F}(\alpha^{\vee}(z+1)))=\Lambda (\tr_{F_{\alpha}/F}(2\overline{a}_{\alpha}z))
\]
for all \(z \in [F_{\alpha}]_{r_{\theta_{\pm}, \alpha}}\).
Note that the factor \(2\) on the right-hand side compensates for the square on the left-hand side, and is harmless since \(p\neq 2\).

For each \(\alpha\in R_{\mathrm{sym}}\) such that
\(r_{\theta_\pm,\alpha}=0\), the extension
\(F_\alpha/F_{\pm\alpha}\) is unramified by {\bf (PD\(_{\alpha}\))}.
Thus, we can choose a trace-zero unit
\(a_\alpha\in F_\alpha^0\cap U_{F_\alpha}\)
compatibly with the Galois action and the relation
\(a_{-\alpha}=-a_{\alpha}\).
Together with lifts, as in \S \ref{sec:modatoas}, of the mod-\(a\)-data
constructed above for the roots of positive depth, these elements
define a genuine function \(a_S\) by \eqref{eq:aSfroma}.

The resulting genuine function is independent of the choices made for the
roots of depth zero. Indeed, if \(a_{\alpha}\) and \(a'_{\alpha}\) are two
trace-zero units, then
\(a'_{\alpha}a_{\alpha}^{-1}\in U_{F_{\pm\alpha}}\),
and \(\kappa_{\alpha}\) is trivial on \(U_{F_{\pm\alpha}}\) since
\(F_\alpha/F_{\pm\alpha}\) is unramified by {\bf (PD\(_{\alpha}\))}.

Strictly speaking, the mod-\(a\)-data for the positive-depth roots,
and hence the resulting genuine function \(a_S\), depend on the choice
of \(\Lambda\), but the character formula in \cite{1907.05173}*{(4.3.2)}
does not.

We remark that although the assumption {\bf (PD)} may not be very explicit in \cite{1907.05173}*{\S4.3}, it is indeed satisfied if the character \(\theta_{\pm}\) comes from a regular supercuspidal parameter as in \cite{1907.05173}*{\S4.3}.

\section{Refined mod-\(a\)-data and canonical splittings over open subgroups} \label{sec:refinedmoda}
We collect a few technical results that will be used in the subsequent construction, especially for \(p=2\).

Let \(\Lambda\) be a non-trivial character of \(F\).
Let \(F_{\pm \alpha}/F\) be a finite separable extension and \(F_{\alpha}/F_{\pm \alpha}\) a quadratic separable extension with the ramification index \(e\).
We put \(d=0\) if \(e=1\) and take \(d\) to be the integer such that the relative discriminant \(\frakd_{F_{\alpha}/F_{\pm \alpha}}\) satisfies \(\frakd_{F_{\alpha}/F_{\pm \alpha}}=\frakp_{F_{\pm \alpha}}^{d+1}\) if \(e=2\).
In particular, we have \(d>0\) if and only if \(F_{\alpha}/F_{\pm \alpha}\) is wildly ramified.

To simplify notation and align with later sections, we use the symbols \(F_{\alpha}/F_{\pm \alpha}\)
to denote a general quadratic extension of a finite separable extension of  \(F\).
We hope this slight abuse of notation will not cause confusion.

Propositions \ref{prop:liftofa}, \ref{prop:wildmoda} below serve as a substitute for \cite{1907.05173}*{Remark 3.5.2} (reviewed in \S \ref{sec:modatoas}) when \(p=2\).

\begin{prop}[cf.\ \cite{1907.05173}*{Proposition D.1.1}] \label{prop:liftofa}
Let $i\leq j$ be integers.
Let \(\xi \colon \frakp_{F_{\alpha}}^i/\frakp_{F_{\alpha}}^{j+1}\to \bbC^{\times}\) be a character.
Let \(r\) be the integer such that \(\Lambda \circ \tr_{F_{\pm \alpha}/F}\) is trivial on \(\frakp_{F_{\pm \alpha}}^{r+1}\) but not on \(\frakp_{F_{\pm \alpha}}^{r}\).
\begin{enumerate}[label=(\arabic*)]
\item \label{prop:liftofa1}
There exists a unique element \(\overline{a} \in \mathfrak{p}_{F_{\alpha}}^{-j+er-d} / \mathfrak{p}_{F_{\alpha}}^{-i+er-d+1}\) such that
we have
\[
\xi (z)=\Lambda (\tr_{F_{\alpha}/F}(\overline{a}z))
\]
for all \(z\in \frakp_{F_{\alpha}}^i\).
\item \label{prop:liftofa2}
Suppose that \(\xi\) is trivial on \(F_{\pm \alpha}\cap \frakp_{F_{\alpha}}^i\).
Then \(\overline{a}\) above lies in the following subset
\[
\frac{\left(F_{\alpha}^0\cap \frakp_{F_{\alpha}}^{-j+er-d}\right)+ \mathfrak{p}_{F_{\alpha}}^{-i+er-d+1}}{\mathfrak{p}_{F_{\alpha}}^{-i+er-d+1}}\subset \mathfrak{p}_{F_{\alpha}}^{-j+er-d} / \mathfrak{p}_{F_{\alpha}}^{-i+er-d+1}.
\]
As this subset is canonically isomorphic to 
\[
\left(F_{\alpha}^0 \cap \frakp_{F_{\alpha}}^{-j+er-d}\right)/\left(F_{\alpha}^0\cap \frakp_{F_{\alpha}}^{-i+er-d+1}\right),
\]
this means that \(\overline{a}\) can be lifted to an element in \(F_{\alpha}^0\cap \frakp_{F_{\alpha}}^{-j+er-d}\).
\end{enumerate}
\end{prop}
\begin{proof}
This is routine;
the proof consists of tracking isomorphisms induced by the auto-duality of \(F_{\alpha}\) and \(F_{\pm \alpha}\) defined by \(\Lambda \circ \tr_{F_{\pm \alpha}/F}\).

\ref{prop:liftofa1} Let \(\left( \mathfrak{p}_{F_{\alpha}}^i / \mathfrak{p}_{F_{\alpha}}^{j+1} \right)^\vee\) be the group of characters of \(\mathfrak{p}_{F_{\alpha}}^i / \mathfrak{p}_{F_{\alpha}}^{j+1}\).
Since \(\Lambda \circ \tr_{F_{\alpha}/F}\) is a character of \(F_{\alpha}\) that is trivial on \(\frakp_{F_{\alpha}}^{er-d+1}\), but not on \(\frakp_{F_{\alpha}}^{er-d}\) (see for example \cite[41.2 Proposition (1), 18.1 Lemma (1)]{MR2234120}),
the map
\(
a \mapsto \left( z \mapsto \Lambda(\mathrm{tr}_{F_{\alpha}/F}(az)) \right)
\)
induces an isomorphism
\[
\mathfrak{p}_{F_{\alpha}}^{-j+er-d} / \mathfrak{p}_{F_{\alpha}}^{-i+er-d+1} \cong \left( \mathfrak{p}_{F_{\alpha}}^i / \mathfrak{p}_{F_{\alpha}}^{j+1} \right)^\vee.
\]
This shows the assertion \ref{prop:liftofa1}.

\ref{prop:liftofa2} We have \(F_{\pm \alpha}\cap \frakp_{F_{\alpha}}^i=\frakp_{F_{\pm \alpha}}^{\lceil \frac{i}{e}\rceil}\).
Via a similar duality for \(F_{\pm \alpha}\), the restriction map is compatible with the trace map:
 \[
 \begin{tikzcd}
\left(\frakp_{F_{\alpha}}^i/\frakp_{F_{\alpha}}^{j+1}\right)^{\vee} 
& \mathfrak{p}_{F_{\alpha}}^{-j+er-d} / \mathfrak{p}_{F_{\alpha}}^{-i+er-d+1}\\
\left(\frakp_{F_{\pm \alpha}}^{\lceil \frac{i}{e}\rceil}/
\frakp_{F_{\pm \alpha}}^{\lceil \frac{j+1}{e}\rceil}\right)^{\vee} 
& \frakp_{F_{\pm \alpha}}^{r+1-\lceil \frac{j+1}{e}\rceil}/\frakp_{F_{\pm \alpha}}^{r+1-\lceil \frac{i}{e}\rceil}. 
\arrow["{\mathrm{res.}}", from=1-1, to=2-1]
\arrow["\simeq"', to=1-1, from=1-2]
\arrow["\operatorname{tr}_{F_{\alpha}/F_{\pm \alpha}}", from=1-2, to=2-2]
\arrow["\simeq"', to=2-1, from=2-2]
 \end{tikzcd}
 \]
Thus, as \(\xi\) is trivial on \(F_{\pm \alpha}\cap \frakp_{F_{\alpha}}^i=\frakp_{F_{\pm \alpha}}^{\lceil \frac{i}{e}\rceil}\), the element \(\overline{a}\) lies in
\[
\left(\frakp_{F_{\alpha}}^{-j+er-d}\cap \tr_{F_{\alpha}/F_{\pm \alpha}}^{-1}\left( \frakp_{F_{\pm \alpha}}^{r+1-\lceil \frac{i}{e}\rceil}\right)\right)/\mathfrak{p}_{F_{\alpha}}^{-i+er-d+1}.
\]
Now, since \(\tr_{F_{\alpha}/F_{\pm \alpha}}(\mathfrak{p}_{F_{\alpha}}^{-i+er-d+1})=\frakp_{F_{\pm \alpha}}^{r+1-\lceil \frac{i}{e}\rceil}\) (see for example \cite[41.2 Proposition (1), 18.1 Lemma (1)]{MR2234120}), we see that
\[
\frakp_{F_{\alpha}}^{-j+er-d}\cap \tr_{F_{\alpha}/F_{\pm \alpha}}^{-1}\left( \frakp_{F_{\pm \alpha}}^{r+1-\lceil \frac{i}{e}\rceil}\right)
=\left(F_{\alpha}^0\cap \frakp_{F_{\alpha}}^{-j+er-d}\right)+ \mathfrak{p}_{F_{\alpha}}^{-i+er-d+1},
\]
and hence the assertion \ref{prop:liftofa2} follows.
\end{proof}
Let \(\tau\) be the generator of the Galois group \(\Gamma_{F_{\alpha}/F_{\pm \alpha}}\) and \(\kappa\) the isomorphism \(F_{\pm \alpha}^{\times}/N_{F_{\alpha}/F_{\pm \alpha}}(F_{\alpha}^{\times}) \stackrel{\sim}{\to} \{\pm 1\}\).
\begin{prop}[cf.\ \cite{1907.05173}*{Proposition D.1.2}] \label{prop:wildmoda}
Let \(m'>m\) be two integers.
Let \(\overline{a} \in (F_{\alpha}^0\cap \frakp_{F_{\alpha}}^{m})/(F_{\alpha}^0\cap \frakp_{F_{\alpha}}^{m'})\) and assume that
some lift \(a\in F_{\alpha}^0\cap \frakp_{F_{\alpha}}^{m}\) (or equivalently, any lift) satisfies \(v_{F_{\alpha}}(a)=m\).

Let \(\delta_{\alpha} \in F_{\alpha}^{\times}\setminus F_{\pm \alpha}^{\times}\) and put \(\delta_{-\alpha}=\tau (\delta_{\alpha})\), so that \(\delta_{\alpha}-\delta_{-\alpha}\) lies in \(F_{\alpha}^0\setminus \{0\}\) and thus for any lift \(a\in F_{\alpha}^0\cap \frakp_{F_{\alpha}}^{m}\) we can evaluate \(\kappa\) at \(\frac{\delta_{\alpha}-\delta_{-\alpha}}{a}\in F_{\pm \alpha}^{\times}\).

Then \(\kappa\left(\frac{\delta_{\alpha}-\delta_{-\alpha}}{a}\right)\) is independent of the choice of a lift \(a\in F_{\alpha}^0\cap \frakp_{F_{\alpha}}^{m}\) of \(\overline{a}\)
if and only if \(m'-m \geq 2d+1\).
\end{prop}
\begin{proof}
If \(a_1, a_2\in F_{\alpha}^0\cap \frakp_{F_{\alpha}}^{m}\) are both lifts of \(\overline{a}\), then we have \(a_2-a_1\in F_{\alpha}^0\cap \frakp_{F_{\alpha}}^{m'}\) and thus
\[
\frac{a_2}{a_1}=1+\frac{a_2-a_1}{a_1}\in U_{F_{\alpha}}^{m'-m}\cap F_{\pm \alpha}^{\times}=U_{F_{\pm \alpha}}^{\lceil \frac{m'-m}{e}\rceil}.
\]
Conversely, for any lift \(a_1\) and \(u\in U_{F_{\pm \alpha}}^{\lceil \frac{m'-m}{e}\rceil}\), the product \(a_1u\) is another lift of \(\overline{a}\).
Hence, the set of all lifts is naturally a torsor under \(U_{F_{\pm \alpha}}^{\lceil \frac{m'-m}{e}\rceil}\).

Now, \(\kappa\left(\frac{\delta_{\alpha}-\delta_{-\alpha}}{a}\right)\) is independent of the choice of a lift \(a\in F_{\alpha}^0\cap \frakp_{F_{\alpha}}^{m}\)
if and only if
\[
\kappa\left(\frac{a_2}{a_1}\right)=1
\]
for any lifts \(a_1, a_2\) of \(\overline{a}\).
This is the case if and only if \(\kappa\) is trivial on \(U_{F_{\pm \alpha}}^{\lceil \frac{m'-m}{e}\rceil}\).
This latter condition is equivalent to \(\lceil \frac{m'-m}{e}\rceil\geq d+1\), or put differently, \(m'-m\geq ed+1\).
Finally, we have \(ed+1=2d+1\), as \(d=0\) unless \(e=2\).
\end{proof}

The following proposition is used in \S \ref{sec:relspl}.
\begin{prop} \label{prop:minramchi}
Suppose that \(d>0\). Let \(\chi\) range over the characters of \(F_{\alpha}^{\times}\) such that \(\chi|_{F_{\pm \alpha}^{\times}}=\kappa\).
Then the minimum of the levels of such \(\chi\) is \(2d\) and \(\chi|_{U_{F_{\alpha}}^{2d}}\) is independent of the choice of \(\chi\) of level \(2d\).
\end{prop}
\begin{proof}
This is easily seen by noting that \(\kappa|_{U_{F_{\pm \alpha}}^d}\) induces a non-trivial character of \(U_{F_{\pm \alpha}}^d/U_{F_{\pm \alpha}}^{d+1}\) and that the inclusion \(U_{F_{\pm \alpha}}^d \hookrightarrow U_{F_{\alpha}}^{2d}\) induces an isomorphism \(U_{F_{\pm \alpha}}^d/U_{F_{\pm \alpha}}^{d+1} \stackrel{\sim}{\to} U_{F_{\alpha}}^{2d}/U_{F_{\alpha}}^{2d+1}\).
\end{proof}
Let us return to the setting in \S \ref{sec:background}.
Given \(\alpha \in R\), we have the associated fields \(F_{\alpha}\), \(F_{\pm \alpha}\) and related objects.
For \(\alpha \in R_{\sym}\), we write \(d_{\alpha}\) for the integer \(d\) associated above to \(F_{\alpha}/F_{\pm \alpha}\). In view of the discussion in \S \ref{sec:dc} on \(\chi\)-data and a splitting of the cover,
for \(\alpha \in R_{\sym}\) with \(d_{\alpha}>0\),
Proposition \ref{prop:minramchi} implies
that the cover \(J_{\alpha}(F)_{\pm}=F_{\alpha}^{\times}/N_{F_{\alpha}/F_{\pm \alpha}}(F_{\alpha}^{\times}) \twoheadrightarrow J_{\alpha}(F)=F_{\alpha}^1\) canonically splits over \(\widetilde{N}_{\alpha}(U_{F_{\alpha}}^{2d})\). 
Similarly, for the other \(\alpha \in R_{\sym}\) with \(d_{\alpha}=0\), the cover canonically splits over \(\widetilde{N}_{\alpha}(U_{F_{\alpha}})\).
(Of course, the cover \(J_{\alpha}(F)_{\pm}\twoheadrightarrow J_{\alpha}(F)\) trivially splits for \(\alpha \in R_{\asym}\).)

\section{Constructions of the genuine function} \label{sec:constr}
Let \(G\) be a connected reductive group and \(S\subset G\) an elliptic maximal torus.
Let \(S(F)_{\pm}\) be the double cover of \(S(F)\) associated to the absolute root system \(R\subset X^{\ast}(S)\) of \(G\) with respect to \(S\) as in \S \ref{sec:dc}.

\subsection{Lift of \(J_{O}(F) \to S(F)\)} \label{sec:liftJOS}
According to \cite{1907.05173}*{Appendix C}, dually to the situation reviewed in \S \ref{sec:dc}, for any \(\Sigma\)-orbit \(O\subset X_{\ast}(S)\), we have a natural homomorphism \(J_{O}\to S\) defined over \(F\),
which we write as \(\mapJOS{O}\).
The following proposition relates \(\mapJOS{O}\) to the composite \(N_{F_{\alpha}/F} \circ \alpha^{\vee}\colon F_{\alpha}^{\times}\to S(F)\) used in \S \ref{sec:gctoas}.
\begin{prop}\label{prop:mapJOS}
	Let \(\alpha \in R\) be a root of \(G\) and \(O\) the \(\Sigma\)-orbit of the coroot \(\alpha^{\vee}\).
	\begin{enumerate}[label=(\arabic*)]
		\item \label{prop:mapJOSasym}
		Suppose that \(\alpha\) is asymmetric. 
		Then the composite \(N_{F_{\alpha}/F} \circ \alpha^{\vee}\colon F_{\alpha}^{\times}\to S(F)\) is equal to
		\[
		\mapJOS{\alpha} \colon F_{\alpha}^{\times}=J_O(F)\to S(F),
		\]
		the homomorphism induced by \(\mapJOS{O}\) and the canonical identification \(F_{\alpha}^{\times}=J_O(F)\).
		\item \label{prop:mapJOSsym}
		Suppose that \(\alpha\) is symmetric. 
		Then the composite \(N_{F_{\alpha}/F} \circ \alpha^{\vee}\colon F_{\alpha}^{\times}\to S(F)\) is equal to the composite of the following two homomorphisms:
		\begin{align*}
			&\widetilde{N}_{\alpha} \colon F_{\alpha}^{\times} \twoheadrightarrow F_{\alpha}^1, \quad x\mapsto \frac{x}{\tau_{\alpha} (x)}, \\
			&\mapJOS{\alpha} \colon F_{\alpha}^1=J_O(F)\to S(F),
		\end{align*}
		where \(\mapJOS{\alpha}\) is the homomorphism induced by \(\mapJOS{O}\) and the canonical identification \(F_{\alpha}^1=J_O(F)\).
		Moreover, we have \(\alpha ^{\vee}\circ \widetilde{N}_{\alpha}=N_{F_{\alpha}/F_{\pm \alpha}} \circ \alpha^{\vee}\colon F_{\alpha}^{\times} \to S(F_{\pm \alpha})\).
	\end{enumerate}
\end{prop}
\begin{proof}
	Left to the reader.
\end{proof}
\begin{rem}
The homomorphism \(\widetilde{N}_{\alpha}\), which already appeared in \S \ref{sec:dc}, has a natural interpretation as the norm map \(N_{F_{\alpha}/F_{\pm \alpha}}\) for the torus \(\Res_{F_{\alpha}/F_{\pm \alpha}}^1 \bbG_m\):
\[
N_{F_{\alpha}/F_{\pm \alpha}} \colon \Res_{F_{\alpha}/F_{\pm \alpha}}^1 \bbG_m(F_{\alpha})=F_{\alpha}^{\times} \twoheadrightarrow \Res_{F_{\alpha}/F_{\pm \alpha}}^1 \bbG_m(F_{\pm \alpha})=F_{\alpha}^1.
\]
\end{rem}

For any \(\Sigma\)-orbit \(O'\subset X^{\ast}(S)\), we similarly write \(\mapSJO{O'} \colon S\to J_{O'}\) for the natural homomorphism used in the construction of \(S(F)_{\pm}\) reviewed in \S \ref{sec:dc}.
\begin{prop}\label{prop:mapSJO}
	Let \(\beta \in R\) be a root of \(G\) and \(O'\) the \(\Sigma\)-orbit of \(\beta\).
	\begin{enumerate}[label=(\arabic*)]
		\item \label{prop:mapSJOasym}
		Suppose that \(\beta\) is asymmetric. 
		We write \(\mapSJO{\beta} \colon S(F) \to F_{\beta}^{\times}\) for the homomorphism induced by \(\mapSJO{O'}\) and the canonical identification \(J_{O'}(F)=F_{\beta}^{\times}\).
		Then \(\mapSJO{\beta}\) is equal to the composite of
		\begin{align*}
			&i_{F_{\beta}/F} \colon S(F) \hookrightarrow S(F_{\beta}), \\
			&\beta \colon S(F_{\beta}) \to F_{\beta}^{\times},
		\end{align*}
		where \(i_{F_{\beta}/F}\) is the inclusion.
		\item \label{prop:mapSJOsym}
		Suppose that \(O'\) is symmetric.
		We write \(\mapSJO{\beta} \colon S(F) \to F_{\beta}^1\) for the homomorphism induced by \(\mapSJO{O'}\) and the canonical identification \(J_{O'}(F)=F_{\beta}^1\).
		Then \(\mapSJO{\beta}\) is equal to the composite of
		\begin{align*}
			&i_{F_{\pm \beta}/F} \colon S(F) \hookrightarrow S(F_{\pm \beta}), \\
			&\beta \colon S(F_{\pm \beta}) \to F_{\beta}^1,
		\end{align*}
		where \(i_{F_{\pm \beta}/F}\) is the inclusion.
		Moreover, we have \(\beta \circ N_{F_{\beta}/F_{\pm \beta}}=\widetilde{N}_{\beta} \circ \beta \colon S(F_{\beta})\to F_{\beta}^1\).
	\end{enumerate}
\end{prop}
\begin{proof}
	Left to the reader.
\end{proof}

\begin{prop}[cf.\ \cite{1907.05173}*{Proposition D.2.1}]
\label{prop:rootwise-lift}
Let \(\alpha\in R\) and let \(O\subset X_{\ast}(S)\) be the \(\Sigma\)-orbit
of \(\alpha^{\vee}\). Then there exists a natural homomorphism
\[
\mapJOS{O, \pm}\colon J_O(F)_\pm\longrightarrow S(F)_\pm
\]
lifting \(\mapJOS{O}\colon J_O(F)\to S(F)\).
When \(O\) is symmetric and \(\alpha \in O\), the homomorphism induced by the canonical identification \(J_{O}(F)_{\pm}=F_{\alpha}^{\times}/N_{F_{\alpha}/F_{\pm \alpha}}(F_{\alpha}^{\times})\)  is the
homomorphism denoted by \(\mapJOS{\alpha, \pm}\) in
\cite{1907.05173}*{Proposition D.2.1}.
\end{prop}
The rest of this subsection is devoted to the proof of Proposition \ref{prop:rootwise-lift}.

We construct a homomorphism \(\mapJOS{O, \pm \pm}\colon J_O(F)_{\pm} \to S(F)_{\pm \pm}\) that lifts \(\mapJOS{O}\), which is then composed with \(S(F)_{\pm \pm} \twoheadrightarrow S(F)_{\pm}\).
By the definition of \(S(F)_{\pm \pm}\), this amounts to constructing a lift \(J_O(F)_{\pm} \to J_{O'}(F)_{\pm}\) of \(\mapSJO{O'} \circ \mapJOS{O} \colon J_O(F)\to J_{O'}(F)\) for each \(\Sigma\)-orbit \(O'\)of a root:
\begin{equation} \label{eq:liftJJ}
	\begin{tikzcd}
		J_O(F)_{\pm} & & J_{O'}(F)_{\pm} \\
		J_O(F) & S(F) & J_{O'}(F).
		\arrow[dashed, from=1-1, to=1-3]
		\arrow[twoheadrightarrow, from=1-1, to=2-1]
		\arrow[twoheadrightarrow, from=1-3, to=2-3]
		\arrow["\mapJOS{O}", from=2-1, to=2-2]
		\arrow["\mapSJO{O'}", from=2-2, to=2-3]
	\end{tikzcd}
\end{equation}
Thus, let \(\beta\) be a root of \(G\) and \(O'\) the \(\Sigma\)-orbit of \(\beta\).

If \(O'\) is asymmetric, we define the lift \(J_O(F)_{\pm} \to J_{O'}(F)_{\pm}\) to be the one given by the canonical section \(J_{O'}(F) \hookrightarrow J_{O'}(F)_{\pm}\).

Suppose that \(O'\) is symmetric.
The main step is to construct a natural homomorphism \(\Nbent{\alpha}{\beta} \colon F_{\alpha}^{\times} \to F_{\beta}^{\times}/N_{F_{\beta}/F_{\pm \beta}}(F_{\beta}^{\times})\) making the following diagram commute:
\begin{equation}\label{eq:Nbent}
	\begin{tikzcd}
		& & & & F_{\beta}^{\times}/N_{F_{\beta}/F_{\pm \beta}}(F_{\beta}^{\times}) \\
		F_{\alpha}^{\times} & S(F_{\alpha}) & S(F) & S(F_{\pm \beta}) & F_{\beta}^1.
		\arrow["\Nbent{\alpha}{\beta}", dashed, bend left=10, from=2-1, to=1-5]
		\arrow["\alpha^{\vee}", near end, from=2-1, to=2-2]
		\arrow["N_{F_{\alpha}/F}", from=2-2, to=2-3]
		\arrow["i_{F_{\pm \beta}/F}", hookrightarrow, from=2-3, to=2-4]
		\arrow["\beta", from=2-4, to=2-5]
		\arrow["\widetilde{N}_{\beta}", twoheadrightarrow, from=1-5, to=2-5]
	\end{tikzcd}
\end{equation}
To rewrite \(\beta \circ i_{F_{\pm \beta}/F} \circ N_{F_{\alpha}/F} \circ \alpha^{\vee}\), we need the following lemma:
\begin{lem} \label{lem:Nbent}
	\begin{enumerate}[label=(\arabic*)]
		\item \label{item:Mackey}
		The composite \(i_{F_{\pm \beta}/F} \circ N_{F_{\alpha}/F} \colon S(F_{\alpha}) \to S(F_{\pm \beta})\) is equal to the product (with respect to the multiplication of \(S(F_{\pm \beta})\)) of the following homomorphisms:
		\[
		N_{gF_{\alpha} \cdot F_{\pm \beta}/F_{\pm \beta}} \circ i_{gF_{\alpha} \cdot F_{\pm \beta}/gF_{\alpha}}\circ g \colon S(F_{\alpha}) \to S(F_{\pm \beta})
		\]
		for \(g \in \Gamma_{\pm \beta} \backslash \Gamma_F/\Gamma_{\alpha}\).
		Here, \(gF_{\alpha} \cdot F_{\pm \beta}\) denotes the composite field.
		\item \label{item:vanish}
		Let \(x\in F_{\alpha}^{\times}\).
        We have
		\[
		\beta \left(N_{gF_{\alpha} \cdot F_{\pm \beta}/F_{\pm \beta}} \left(g (\alpha^{\vee}(x))\right)\right)=1
		\]
		for any \(g\in \Gamma_F\) such that \(gF_{\alpha} \cdot F_{\pm \beta}\) does not contain \(F_{\beta}\).
	\end{enumerate}
\end{lem}
\begin{proof}
	The assertion \ref{item:Mackey} is an analogue of the Mackey decomposition.
	Indeed, let \(\gamma \in S(F_{\alpha})\). Rearranging the product, we get
\begin{align*}
i_{F_{\pm \beta}/F} \circ N_{F_{\alpha}/F}(\gamma)
=\prod_{\sigma \in \Gamma_F/\Gamma_{\alpha}} \sigma (\gamma)
&=\prod_{g \in \Gamma_{\pm \beta} \backslash \Gamma_F/\Gamma_{\alpha}} \prod_{h\in \Gamma_{\pm \beta}/\left( \Gamma_{g\alpha} \cap \Gamma_{\pm \beta}\right)} hg (\gamma) \\
&=\prod_{g \in \Gamma_{\pm \beta} \backslash \Gamma_F/\Gamma_{\alpha}} N_{gF_{\alpha} \cdot F_{\pm \beta}/F_{\pm \beta}} \left(g (\gamma)\right).
\end{align*}
	
	For \(g\in \Gamma_F\) as in \ref{item:vanish},
	there exists \(\tau \in (\Gamma_{g\alpha}\cap \Gamma_{\pm \beta}) \setminus \Gamma_{\beta}\).
	On the other hand, we have
	\begin{align*}
		\beta \left(N_{gF_{\alpha} \cdot F_{\pm \beta}/F_{\pm \beta}} \left(g (\alpha^{\vee}(x))\right)\right)
		&=\prod_{h\in \Gamma_{\pm \beta}/\left( \Gamma_{g\alpha} \cap \Gamma_{\pm \beta}\right)} \beta \left(h\left(g (\alpha^{\vee}(x))\right)\right) \\
		&=\prod_{h\in \Gamma_{\pm \beta}/\left( \Gamma_{g\alpha} \cap \Gamma_{\pm \beta}\right)} h\left({}^{h^{-1}}\beta \circ {}^{g}\alpha^{\vee}(g(x))\right),
	\end{align*}
	where we write \({}^{g}\alpha ^{\vee}=g\alpha ^{\vee}=g\circ \alpha ^{\vee} \circ g^{-1}\) and similarly \({}^{h^{^{-1}}}\beta=h^{-1}\beta=h^{-1}\circ \beta \circ h\) for clarity.
	However, we have \({}^{h^{-1}}\beta \circ {}^{g}\alpha^{\vee}=0\) in \(X^{\ast}(\bbG_m)=\bbZ\), as
	\begin{align*}
		{}^{h^{-1}}\beta \circ {}^{g}\alpha^{\vee}
		&=\langle {}^{h^{-1}}\beta, {}^{g}\alpha^{\vee}\rangle \\
		&=\langle {}^{h^{-1}}\beta, {}^{\tau g}\alpha^{\vee}\rangle \\
		&=\langle {}^{\tau ^{-1}h^{-1}}\beta, {}^{g}\alpha^{\vee}\rangle \\
		&=\langle -{}^{h^{-1}}\beta, {}^{g}\alpha^{\vee}\rangle 
		=-{}^{h^{-1}}\beta \circ {}^{g}\alpha^{\vee}.
	\end{align*}
\end{proof}
By Lemma \ref{lem:Nbent}, we see that
\begin{align*}
	\beta \circ i_{F_{\pm \beta}/F} \circ N_{F_{\alpha}/F} \circ \alpha^{\vee}(x)
	&=\prod_{g \in \Gamma_{\pm \beta} \backslash \Gamma_F/\Gamma_{\alpha}} \beta \left(N_{gF_{\alpha} \cdot F_{\pm \beta}/F_{\pm \beta}} \left(g (\alpha^{\vee}(x))\right)\right) \\
	&=\prod_{\substack{g \in \Gamma_{\pm \beta} \backslash \Gamma_F/\Gamma_{\alpha} \\ gF_{\alpha} \cdot F_{\pm \beta} \supset F_{\beta}}} \beta \left(N_{gF_{\alpha} \cdot F_{\pm \beta}/F_{\pm \beta}} \left(g (\alpha^{\vee}(x))\right)\right) \\
	&=\prod_{\substack{g \in \Gamma_{\pm \beta} \backslash \Gamma_F/\Gamma_{\alpha} \\ gF_{\alpha} \cdot F_{\pm \beta} \supset F_{\beta}}} \beta \left(N_{F_{\beta}/F_{\pm \beta}}\circ N_{gF_{\alpha} \cdot F_{\pm \beta}/F_{\beta}} \left(g (\alpha^{\vee}(x))\right)\right) \\
    &=\prod_{\substack{g \in \Gamma_{\pm \beta} \backslash \Gamma_F/\Gamma_{\alpha} \\ gF_{\alpha} \cdot F_{\pm \beta} \supset F_{\beta}}} \widetilde{N}_{\beta} \left( \beta \left( N_{gF_{\alpha} \cdot F_{\pm \beta}/F_{\beta}} \left(g (\alpha^{\vee}(x))\right)\right) \right)\\
    &=\widetilde{N}_{\beta}\left(\prod_{\substack{g \in \Gamma_{\pm \beta} \backslash \Gamma_F/\Gamma_{\alpha} \\ gF_{\alpha} \cdot F_{\pm \beta} \supset F_{\beta}}} \beta \left( N_{gF_{\alpha} \cdot F_{\pm \beta}/F_{\beta}} \left(g (\alpha^{\vee}(x))\right)\right)\right).
\end{align*}
It is now clear that the following homomorphism gives \(\Nbent{\alpha}{\beta}\) in the diagram \eqref{eq:Nbent}:
\begin{equation} \label{eq:defNbent}
	\Nbent{\alpha}{\beta} \colon F_{\alpha}^{\times} \to F_{\beta}^{\times}/N_{F_{\beta}/F_{\pm \beta}}(F_{\beta}^{\times}), \quad 
	x\mapsto
	\prod_{\substack{g \in \Gamma_{\pm \beta} \backslash \Gamma_F/\Gamma_{\alpha} \\ gF_{\alpha} \cdot F_{\pm \beta} \supset F_{\beta}}} \beta \left( N_{gF_{\alpha} \cdot F_{\pm \beta}/F_{\beta}} \left(g (\alpha^{\vee}(x))\right)\right).
\end{equation}

Suppose that \(O\) is asymmetric. Then, under the canonical identifications \(J_{O}(F)_{\pm}=F_{\alpha}^{\times} \times \{\pm 1\}, J_{O}(F)=F_{\alpha}^{\times}\), \(J_{O'}(F)_{\pm}=F_{\beta}^{\times}/N_{F_{\beta}/F_{\pm \beta}}(F_{\beta}^{\times}), J_{O'}(F)=F_{\beta}^1\), the diagram \eqref{eq:liftJJ} reads
\[
\begin{tikzcd}
	F_{\alpha}^{\times} \times \{\pm 1\} & & F_{\beta}^{\times}/N_{F_{\beta}/F_{\pm \beta}}(F_{\beta}^{\times}) \\
	F_{\alpha}^{\times} & S(F) & F_{\beta}^1.
	\arrow[dashed, from=1-1, to=1-3]
	\arrow[twoheadrightarrow, from=1-1, to=2-1]
	\arrow[twoheadrightarrow, from=1-3, to=2-3]
	\arrow["\mapJOS{\alpha}", from=2-1, to=2-2]
	\arrow["\mapSJO{\beta}", from=2-2, to=2-3]
\end{tikzcd}
\]
Since \(\mapSJO{\beta} \circ \mapJOS{\alpha}=\beta \circ i_{F_{\pm \beta}/F} \circ N_{F_{\alpha}/F} \circ \alpha ^{\vee}\) by Propositions \ref{prop:mapJOS} \ref{prop:mapJOSasym} and \ref{prop:mapSJO} \ref{prop:mapSJOsym}, we compose \(F_{\alpha}^{\times} \times \{\pm 1\} \twoheadrightarrow F_{\alpha}^{\times}\) with \(\Nbent{\alpha}{\beta}\) to get the dashed arrow.
It is routine to check that the induced homomorphism \(J_O(F)_{\pm}\to J_{O'}(F)_{\pm}\) is independent of the choice of \(\alpha, \beta\).

Suppose that \(O\) is symmetric.
This time \(J_{O}(F)_{\pm}\) and \(J_{O}(F)\) are identified with \(F_{\alpha}^{\times}/N_{F_{\alpha}/F_{\pm \alpha}}(F_{\alpha}^{\times})\) and \(F_{\alpha}^1\) respectively, and the diagram \eqref{eq:liftJJ} reads
\begin{equation} \label{eq:liftJJsym}
\begin{tikzcd}
	F_{\alpha}^{\times}/N_{F_{\alpha}/F_{\pm \alpha}}(F_{\alpha}^{\times}) & & F_{\beta}^{\times}/N_{F_{\beta}/F_{\pm \beta}}(F_{\beta}^{\times}) \\
	F_{\alpha}^1 & S(F) & F_{\beta}^1.
	\arrow[dashed, from=1-1, to=1-3]
	\arrow["\widetilde{N}_{\alpha}"', twoheadrightarrow, from=1-1, to=2-1]
	\arrow["\widetilde{N}_{\beta}", twoheadrightarrow, from=1-3, to=2-3]
	\arrow["\mapJOS{\alpha}", from=2-1, to=2-2]
	\arrow["\mapSJO{\beta}", from=2-2, to=2-3]
\end{tikzcd}
\end{equation}
Since \(\mapSJO{\beta} \circ \mapJOS{\alpha} \circ \widetilde{N}_{\alpha}=\beta \circ i_{F_{\pm \beta}/F} \circ N_{F_{\alpha}/F} \circ \alpha ^{\vee}\) by Propositions \ref{prop:mapJOS} \ref{prop:mapJOSsym} and \ref{prop:mapSJO} \ref{prop:mapSJOsym},
it suffices to show the following lemma:
\begin{lem} \label{lem:Nbentfactor}
Let \(\alpha, \beta \in R\) be symmetric as above.
Then, the homomorphism \(\Nbent{\alpha}{\beta}\) in \eqref{eq:Nbent} factors through \(F_{\alpha}^{\times} \twoheadrightarrow F_{\alpha}^{\times}/N_{F_{\alpha}/F_{\pm \alpha}}(F_{\alpha}^{\times})\).
\end{lem}
\begin{proof}
Let \(x\in N_{F_{\alpha}/F_{\pm \alpha}}(F_{\alpha}^{\times})\) and write \(x=N_{F_{\alpha}/F_{\pm \alpha}}(y)\) for some \(y\in F_{\alpha}^{\times}\).
We fix \(\tau_{\alpha} \in \Gamma_{\pm \alpha} \setminus \Gamma_{\alpha}\) and compute the product in the definition \eqref{eq:defNbent} of \(\Nbent{\alpha}{\beta}(x)\) by splitting the index set into three parts according to the following conditions:
\begin{itemize}
\item \(\Gamma_{\pm \beta} g\Gamma_{\alpha} \neq \Gamma_{\pm \beta} g\tau_{\alpha} \Gamma_{\alpha}\);
\item \(\Gamma_{\pm \beta} g\Gamma_{\alpha}=\Gamma_{\pm \beta} g\tau_{\alpha} \Gamma_{\alpha}\) and \(\Gamma_{\beta} g\Gamma_{\alpha}=\Gamma_{\beta}g\tau_{\alpha}\Gamma_{\alpha}\);
\item \(\Gamma_{\pm \beta} g\Gamma_{\alpha}=\Gamma_{\pm \beta} g\tau_{\alpha} \Gamma_{\alpha}\) and \(\Gamma_{\beta} g\Gamma_{\alpha}\neq \Gamma_{\beta}g\tau_{\alpha}\Gamma_{\alpha}\). 
\end{itemize}

First, let \(g\in \Gamma_F\) satisfy \(\Gamma_{\pm \beta} g\Gamma_{\alpha} \neq \Gamma_{\pm \beta} g\tau_{\alpha} \Gamma_{\alpha}\)
as well as \(gF_{\alpha} \cdot F_{\pm \beta} \supset F_{\beta}\).
Since \(g\tau_{\alpha}F_{\alpha}=gF_{\alpha}\), we also have \(g\tau_{\alpha}F_{\alpha} \cdot F_{\pm \beta} \supset F_{\beta}\).
As \(x\in F_{\pm \alpha}^{\times}\), we have \(\tau_{\alpha}(\alpha^{\vee}(x))=\alpha^{\vee}(x)^{-1}\) and hence
\begin{align*}
&\beta \left( N_{gF_{\alpha} \cdot F_{\pm \beta}/F_{\beta}} \left(g (\alpha^{\vee}(x))\right)\right) \cdot \beta \left( N_{g\tau_{\alpha}F_{\alpha} \cdot F_{\pm \beta}/F_{\beta}} \left(g\tau_{\alpha} (\alpha^{\vee}(x))\right)\right) \\
=& \beta \left( N_{gF_{\alpha} \cdot F_{\pm \beta}/F_{\beta}} \left(g (\alpha^{\vee}(x))\right)\right) \cdot \beta \left( N_{gF_{\alpha} \cdot F_{\pm \beta}/F_{\beta}} \left(g (\alpha^{\vee}(x))\right)\right)^{-1}=1.
\end{align*}
Thus, by pairing the terms for \(g\) and \(g\tau_{\alpha}\), we see that
\begin{equation} \label{eq:prod1}
\prod_{\substack{g \in \Gamma_{\pm \beta} \backslash \Gamma_F/\Gamma_{\alpha} \\ gF_{\alpha} \cdot F_{\pm \beta} \supset F_{\beta}
\\ \Gamma_{\pm \beta} g\Gamma_{\alpha} \neq \Gamma_{\pm \beta} g\tau_{\alpha} \Gamma_{\alpha}}} \beta \left( N_{gF_{\alpha} \cdot F_{\pm \beta}/F_{\beta}} \left(g (\alpha^{\vee}(x))\right)\right)=1 \in F_{\beta}^{\times}.
\end{equation}

Before treating the other two cases, we compute
\begin{equation} \label{eq:norm}
\begin{split}
g(\alpha^{\vee}(x))
=g(\alpha^{\vee}(N_{F_{\alpha}/F_{\pm \alpha}}(y))) 
=&g(\alpha^{\vee}(y)) \cdot g(\alpha^{\vee}(\tau_{\alpha}(y))) \\
=&g(\alpha^{\vee}(y)) \cdot g\tau_{\alpha}(\alpha^{\vee}(y))^{-1}
\end{split}
\end{equation}
for any \(g\in \Gamma_F\).

Next, let \(g\in \Gamma_F\) satisfy \(\Gamma_{\pm \beta} g\Gamma_{\alpha}=\Gamma_{\pm \beta} g\tau_{\alpha} \Gamma_{\alpha}\) and \(\Gamma_{\beta} g\Gamma_{\alpha}=\Gamma_{\beta}g\tau_{\alpha}\Gamma_{\alpha}\) (note that the former is implied by the latter)
as well as \(gF_{\alpha} \cdot F_{\pm \beta} \supset F_{\beta}\).
Then, we have \(g\tau_{\alpha}=\gamma_1 g \gamma_2\) for some \(\gamma_1 \in \Gamma_{\beta}\) and \(\gamma_2 \in \Gamma_{\alpha}\).
Using this, we compute
\begin{align*}
N_{gF_{\alpha} \cdot F_{\pm \beta}/F_{\beta}} \left(g\tau_{\alpha} (\alpha^{\vee}(y))\right) 
=& N_{gF_{\alpha} \cdot F_{\pm \beta}/F_{\beta}} \left(\gamma_1 g \gamma_2(\alpha^{\vee}(y))\right) \\
=& N_{gF_{\alpha} \cdot F_{\pm \beta}/F_{\beta}} \left(g (\alpha^{\vee}(y))\right).
\end{align*}
Thus, we have
\(N_{gF_{\alpha}\cdot F_{\pm \beta}/F_{\beta}}(g(\alpha^{\vee}(x))) 
=1\)
by \eqref{eq:norm},
and so
\begin{equation} \label{eq:prod2}
\prod_{\substack{g \in \Gamma_{\pm \beta} \backslash \Gamma_F/\Gamma_{\alpha} \\ gF_{\alpha} \cdot F_{\pm \beta} \supset F_{\beta}
\\ \Gamma_{\pm \beta} g\Gamma_{\alpha}=\Gamma_{\pm \beta} g\tau_{\alpha} \Gamma_{\alpha}
\\ \Gamma_{\beta} g\Gamma_{\alpha}=\Gamma_{\beta}g\tau_{\alpha}\Gamma_{\alpha}}} \beta \left( N_{gF_{\alpha} \cdot F_{\pm \beta}/F_{\beta}} \left(g (\alpha^{\vee}(x))\right)\right)=1 \in F_{\beta}^{\times}.
\end{equation}

Finally, let \(g\in \Gamma_F\) satisfy \(\Gamma_{\pm \beta} g\Gamma_{\alpha}=\Gamma_{\pm \beta} g\tau_{\alpha} \Gamma_{\alpha}\) and \(\Gamma_{\beta} g\Gamma_{\alpha}\neq \Gamma_{\beta}g\tau_{\alpha}\Gamma_{\alpha}\)
as well as \(gF_{\alpha} \cdot F_{\pm \beta} \supset F_{\beta}\).
Then, we have \(g\tau_{\alpha}=\tau_{\beta}g\gamma\) for some \(\tau_{\beta} \in \Gamma_{\pm \beta} \setminus \Gamma_{\beta}\) and \(\gamma \in \Gamma_{\alpha}\).
Using this, we compute
\[
g\tau_{\alpha} (\alpha^{\vee}(y))
=\tau_{\beta} g \gamma(\alpha^{\vee}(y))
=\tau_{\beta} g (\alpha^{\vee}(y)).
\]
Abbreviating \(z=g(\alpha^{\vee}(y)) \in S(gF_{\alpha})\),
we combine the above with \eqref{eq:norm} to find that
\begin{align*}
\beta \left(N_{gF_{\alpha}\cdot F_{\pm \beta}/F_{\beta}}(g(\alpha^{\vee}(x)))\right) 
=& \beta \left(N_{gF_{\alpha}\cdot F_{\pm \beta}/F_{\beta}}(z)\right) \cdot \beta \left(N_{gF_{\alpha}\cdot F_{\pm \beta}/F_{\beta}}(\tau_{\beta} (z))\right)^{-1} \\
=& \beta \left(\prod_{h\in \Gamma_{\beta}/\left(\Gamma_{g\alpha} \cap \Gamma_{\pm \beta}\right)}h(z)\right) \cdot \beta \left(\prod_{h\in \Gamma_{\beta}/\left(\Gamma_{g\alpha} \cap \Gamma_{\pm \beta}\right)}h(\tau_{\beta} (z))\right)^{-1} \\
=& \prod_{h\in \Gamma_{\beta}/(\Gamma_{g\alpha} \cap \Gamma_{\pm \beta})}h(\beta (z)) \cdot \prod_{h\in \Gamma_{\beta}/\left(\Gamma_{g\alpha} \cap \Gamma_{\pm \beta}\right)}h(\beta (\tau_{\beta} (z)))^{-1} \\
=& \prod_{h\in \Gamma_{\beta}/(\Gamma_{g\alpha} \cap \Gamma_{\pm \beta})}h(\beta (z)) \cdot \prod_{h\in \Gamma_{\beta}/\left(\Gamma_{g\alpha} \cap \Gamma_{\pm \beta}\right)}h\tau_{\beta}(\beta (z)) \\
=& \prod_{h\in \Gamma_{\pm \beta}/(\Gamma_{g\alpha} \cap \Gamma_{\pm \beta})}h(\beta (z)) \\
=& N_{gF_{\alpha}\cdot F_{\pm \beta}/F_{\pm \beta}}(\beta (z)) 
\in N_{F_{\beta}/F_{\pm \beta}}(F_{\beta}^{\times}),
\end{align*}
and thus
\begin{equation} \label{eq:prod3}
\prod_{\substack{g \in \Gamma_{\pm \beta} \backslash \Gamma_F/\Gamma_{\alpha} \\ gF_{\alpha} \cdot F_{\pm \beta} \supset F_{\beta} \\ \Gamma_{\pm \beta} g\Gamma_{\alpha}=\Gamma_{\pm \beta} g\tau_{\alpha} \Gamma_{\alpha} \\ \Gamma_{\beta} g\Gamma_{\alpha}\neq \Gamma_{\beta}g\tau_{\alpha}\Gamma_{\alpha}}} \beta \left( N_{gF_{\alpha} \cdot F_{\pm \beta}/F_{\beta}} \left(g (\alpha^{\vee}(x))\right)\right)
\in N_{F_{\beta}/F_{\pm \beta}}(F_{\beta}^{\times}).
\end{equation}

Now, by \eqref{eq:prod1}, \eqref{eq:prod2}, and \eqref{eq:prod3}, we conclude that \(\Nbent{\alpha}{\beta}(x)\) is indeed trivial in \(F_{\beta}^{\times}/N_{F_{\beta}/F_{\pm \beta}}(F_{\beta}^{\times})\).
\end{proof}
By Lemma \ref{lem:Nbentfactor}, we see that the homomorphism \(\Nbent{\alpha}{\beta}\) (with \(F_{\alpha}^{\times}/N_{F_{\alpha}/F_{\pm \alpha}}(F_{\alpha}^{\times})\) viewed as the domain)
serves as the dashed arrow in \eqref{eq:liftJJsym}.
Again it is routine to check that the induced homomorphism \(J_O(F)_{\pm}\to J_{O'}(F)_{\pm}\) is independent of the choice of \(\alpha, \beta\).

We have thus constructed a lift \(J_{O}(F)_{\pm} \to J_{O'}(F)_{\pm}\) of \(\mapSJO{O'} \circ \mapJOS{O} \colon J_{O}(F)\to J_{O'}(F)\) for each \(\Sigma\)-orbit \(O'\). These yield a homomorphism \(\mapJOS{O, \pm \pm } \colon J_{O}(F)_{\pm} \to S(F)_{\pm \pm}\), which in turn is composed with \(S(F)_{\pm \pm} \twoheadrightarrow S(F)_{\pm}\) to give the desired lift \(\mapJOS{O, \pm} \colon J_{O}(F)_{\pm} \to S(F)_{\pm}\) of \(\mapJOS{O} \colon J_{O}(F) \to S(F)\).
This completes the proof of Proposition \ref{prop:rootwise-lift}.
\qed

\subsection{Method 1: Direct construction without exponential map} \label{sec:method1}
Let \(\theta_{\pm}\) be a genuine character of \(S(F)_{\pm}\) and \(\Lambda\) be a non-trivial additive character of \(F\).

For \(\alpha \in R_{\sym}\), 
we write \(\mapJOS{\alpha, \pm} \colon F_{\alpha}^{\times}/N_{F_{\alpha}/F_{\pm \alpha}}(F_{\alpha}^{\times}) \to S(F)_{\pm}\)
for the homomorphism induced by \(\mapJOS{O, \pm}\) and the canonical identification \(J_{O}(F)_{\pm}=F_{\alpha}^{\times}/N_{F_{\alpha}/F_{\pm \alpha}}(F_{\alpha}^{\times})\).
We write \(\theta_{\pm, \alpha}\) for the composite of \(\theta_{\pm} \circ \mapJOS{\alpha, \pm}\) with the canonical surjection \(F_{\alpha}^{\times} \twoheadrightarrow F_{\alpha}^{\times}/N_{F_{\alpha}/F_{\pm \alpha}}(F_{\alpha}^{\times})\).

For \(\alpha \in R_{\sym}\), let \(e_{\alpha}\) be the ramification index of \(F_{\alpha}/F_{\pm \alpha}\).
Let \(d_{\alpha}\) be the integer such that the relative discriminant \(\frakd_{F_{\alpha}/F_{\pm \alpha}}\) equals \(\frakp_{F_{\pm \alpha}}^{d_{\alpha}+1}\) if \(e_{\alpha}=2\) and put \(d_{\alpha}=0\) if \(e_{\alpha}=1\).
Let \(r_{\alpha}\) be the integer such that \(\Lambda \circ \tr_{F_{\pm \alpha}/F}\) is trivial on \(\frakp_{F_{\pm \alpha}}^{r_{\alpha}+1}\) but non-trivial on \(\frakp_{F_{\pm \alpha}}^{r_{\alpha}}\).

For \(\alpha \in R_{\sym}\), 
let \(l_{\alpha}\) be the level of \(\theta_{\pm, \alpha}\).
We consider the following condition:
\begin{description}
	\item[(Ram\(_{\alpha}\))]
	if \(F_{\alpha}/F_{\pm \alpha}\) is ramified, then \(l_{\alpha}\geq 4d_{\alpha}+1\).	
\end{description}
Note that if \(p\) is odd and hence \(d_{\alpha}=0\), then {\bf (Ram\(_{\alpha}\))} is equivalent to {\bf (PD\(_{\alpha}\))} in \S \ref{sec:gctoas};
this can be proved in the same way as Theorem \ref{thm:main1} \ref{item:coincidence1} below.
We write {\bf (Ram)} for the condition that {\bf (Ram\(_{\alpha}\))} holds for every \(\alpha \in R_{\sym}\).
\begin{rem} \label{rem:Ramandast}
Recall the condition
\begin{description}
\item[\((\ast)\)] if \(\alpha\) is symmetric and \(d_{\alpha}>0\), then \(l_{\alpha}\geq 4d_{\alpha}+1\)
\end{description}
in \cite{1907.05173}*{Theorem D.3.1}.
The condition {\bf (Ram)} is slightly stronger than the condition \((\ast)\).
When \(p=2\), the two conditions are equivalent, since a ramified quadratic extension is then wildly ramified. When \(p\neq2\), however, \(d_{\alpha}=0\), so that \((\ast)\) imposes no condition, whereas {\bf (Ram\(_{\alpha}\))} requires \(l_\alpha>0\) whenever \(F_{\alpha}/F_{\pm\alpha}\) is ramified.
As noted above, this is precisely the positive-depth condition {\bf (PD\(_{\alpha}\))} of \S \ref{sec:gctoas}. The need for this additional condition in odd residual characteristic was overlooked in \cite{1907.05173}*{Appendix D}.
\end{rem}

\begin{thm}[cf.\ \cite{1907.05173}*{Theorem D.3.1}] \label{thm:main1}
Let \(\theta_{\pm}\) be a genuine character and assume that it satisfies {\bf (Ram)}.

Then the following hold:
\begin{enumerate}[label=(\arabic*)]
\item \label{item:linearization1}
For \(\alpha \in R_{\sym}\) with \(l_{\alpha} >0\),
there exists a unique element
\[
\overline{a}_{\alpha} \in \frakp_{F_{\alpha}}^{-l_{\alpha}+e_{\alpha}r_{\alpha}-d_{\alpha}}/\frakp_{F_{\alpha}}^{-l_{\alpha}+e_{\alpha}r_{\alpha}+d_{\alpha}+1}
\]
such that
\begin{equation} \label{eq:linearization4}
	\theta_{\pm, \alpha}(1+z)=\Lambda (\tr_{F_{\alpha}/F}(\overline{a}_{\alpha}z))
\end{equation}
holds for any \(z\in \frakp_{F_{\alpha}}^{l_{\alpha}-2d_{\alpha}}\).
\item \label{item:liftofa1}
For \(\alpha \in R_{\sym}\) with \(l_{\alpha} >0\),
the element \(\overline{a}_{\alpha}\) defined by \eqref{eq:linearization4} can be lifted to an element \(a_{\alpha} \in F_{\alpha}^0\cap \frakp_{F_{\alpha}}^{-l_{\alpha}+e_{\alpha}r_{\alpha}-d_{\alpha}}\).
For \(\alpha \in R_{\sym}\) with \(l_{\alpha}=0\), the extension \(F_{\alpha}/F_{\pm \alpha}\) is unramified by {\bf (Ram\(_{\alpha}\))}
and we choose a trace-zero unit \(a_{\alpha} \in F_{\alpha}^0 \cap U_{F_{\alpha}}\).
These choices define a genuine function \(a_S\) by \eqref{eq:aSfroma} as in \S \ref{sec:gctoas}.
\item \label{item:independence1}
The genuine function \(a_S\) is independent of all the choices in \ref{item:liftofa1}.
\item \label{item:coincidence1}
If \(p\neq 2\), then the genuine function defined this way coincides with the one defined by the mod-\(a\)-data of \cite{1907.05173}*{\S 4.3} recalled in \S \ref{sec:gctoas}.
\end{enumerate}
\end{thm}
\begin{proof}
For \(\alpha \in R_{\sym}\) with \(l_{\alpha} >0\), we put
\[
i_{\alpha}:=l_{\alpha}-2d_{\alpha}, \quad j_{\alpha}:=l_{\alpha},
\]
so that \(-j_{\alpha}+e_{\alpha}r_{\alpha}-d_{\alpha}=-l_{\alpha}+e_{\alpha}r_{\alpha}-d_{\alpha}\) and
\(-i_{\alpha}+e_{\alpha}r_{\alpha}-d_{\alpha}+1=-l_{\alpha}+e_{\alpha}r_{\alpha}+d_{\alpha}+1\).

To prove the assertions \ref{item:linearization1}, \ref{item:liftofa1}, take \(\alpha \in R_{\sym}\) with \(l_{\alpha} >0\).
We note that \(i_{\alpha}\) is positive
and that we have a usual isomorphism \(\frakp_{F_{\alpha}}^{i_{\alpha}}/\frakp_{F_{\alpha}}^{j_{\alpha}+1} \simeq U_{F_{\alpha}}^{i_{\alpha}}/U_{F_{\alpha}}^{j_{\alpha}+1}; z\mapsto 1+z\).
This is clear if \(d_{\alpha}=0\).
For \(\alpha\) with \(d_{\alpha}>0\), we check that
\[
2i_{\alpha}-(j_{\alpha}+1)=2(l_{\alpha}-2d_{\alpha})-(l_{\alpha}+1)=l_{\alpha}-4d_{\alpha}-1\geq 0
\]
by {\bf (Ram\(_{\alpha}\))}.
Thus, by composing it with \(\theta_{\pm, \alpha}|_{U_{F_{\alpha}}^{i_{\alpha}}}\), we get a character of  \(\frakp_{F_{\alpha}}^{i_{\alpha}}/\frakp_{F_{\alpha}}^{j_{\alpha}+1}\).
Now Proposition \ref{prop:liftofa} \ref{prop:liftofa1} implies the assertion \ref{item:linearization1}.

To prove the assertion \ref{item:liftofa1}, 
we claim that \(U_{F_{\alpha}}^{i_{\alpha}}\cap F_{\pm \alpha}^{\times}\) is contained in \(N_{F_{\alpha}/F_{\pm \alpha}}(F_{\alpha}^{\times})\), and thus \(\theta_{\pm, \alpha}\) is trivial on \(U_{F_{\alpha}}^{i_{\alpha}}\cap F_{\pm \alpha}^{\times}\).
If \(F_{\alpha}/F_{\pm \alpha}\) is unramified, then indeed
\(U_{F_{\alpha}}^{i_{\alpha}}\cap F_{\pm \alpha}^{\times}=U_{F_{\pm \alpha}}^{i_{\alpha}} \subset U_{F_{\pm \alpha}} \subset N_{F_{\alpha}/F_{\pm \alpha}}(F_{\alpha}^{\times})\).
If \(F_{\alpha}/F_{\pm \alpha}\) is ramified, then 
we check that
\[
\left\lceil \frac{i_{\alpha}}{2}\right\rceil-(d_{\alpha}+1)
=\left\lceil \frac{i_{\alpha}-2d_{\alpha}}{2}\right\rceil-1
=\left\lceil \frac{l_{\alpha}-4d_{\alpha}}{2}\right\rceil-1
\geq 0
\]
by {\bf (Ram\(_{\alpha}\))} to find that
\[
U_{F_{\alpha}}^{i_{\alpha}}\cap F_{\pm \alpha}^{\times}=U_{F_{\pm \alpha}}^{\lceil \frac{i_{\alpha}}{2}\rceil} \subset 
U_{F_{\pm \alpha}}^{d_{\alpha}+1} \subset N_{F_{\alpha}/F_{\pm \alpha}}(F_{\alpha}^{\times}),
\]
as claimed.
Then Proposition \ref{prop:liftofa} \ref{prop:liftofa2} implies that \(\overline{a}_{\alpha}\) can indeed be lifted to an element \(a_{\alpha} \in F_{\alpha}^0\cap \frakp_{F_{\alpha}}^{-j_{\alpha}+e_{\alpha}r_{\alpha}-d_{\alpha}}\).
If \(l_{\alpha}=0\), then \(F_{\alpha}/F_{\pm \alpha}\) is unramified by {\bf (Ram\(_{\alpha}\))}, and hence \(F_{\alpha}^0\) contains a unit.
Then we obtain a genuine function \(a_S\) by the formula \eqref{eq:aSfroma}.

Let us prove the assertion \ref{item:independence1}.
For \(\alpha \in R_{\sym}\) such that \(l_{\alpha}>0\), we have \(v_{F_{\alpha}}(a_{\alpha})=-j_{\alpha}+e_{\alpha}r_{\alpha}-d_{\alpha}\), since the level of \(\theta_{\pm, \alpha} \) is \(l_{\alpha}=j_{\alpha}\).
Proposition \ref{prop:wildmoda} shows that \(\kappa_{\alpha} \left(\frac{\delta_{\alpha}-\delta_{-\alpha}}{a_{\alpha}}\right)\) in the formula \eqref{eq:aSfroma} is independent of the choice of the lift \(a_{\alpha}\in F_{\alpha}^0\cap \frakp_{F_{\alpha}}^{-j_{\alpha}+e_{\alpha}r_{\alpha}-d_{\alpha}}\)
of \(\overline{a}_{\alpha} \in \frakp_{F_{\alpha}}^{-j_{\alpha}+e_{\alpha}r_{\alpha}-d_{\alpha}}/\frakp_{F_{\alpha}}^{-i_{\alpha}+e_{\alpha}r_{\alpha}-d_{\alpha}+1}\), as
\[
(-i_{\alpha}+e_{\alpha}r_{\alpha}-d_{\alpha}+1)-(-j_{\alpha}+e_{\alpha}r_{\alpha}-d_{\alpha}) 
=j_{\alpha}-i_{\alpha}+1=l_{\alpha}-(l_{\alpha}-2d_{\alpha})+1=2d_{\alpha}+1.
\]
Moreover, for \(\alpha \in R_{\sym}\) such that \(l_{\alpha}=0\), 
any two choices \(a_{\alpha}, a_{\alpha}'\) satisfy \(a_{\alpha}'/a_{\alpha} \in U_{F_{\pm \alpha}}\), on which \(\kappa_{\alpha}\) is trivial since \(F_{\alpha}/F_{\pm \alpha}\) is unramified.
Thus, we have proved the assertion \ref{item:independence1}.

Assume that \(p\neq 2\) and let us finally prove the assertion \ref{item:coincidence1}.
Recall that the mod-\(a\)-data used in \S \ref{sec:gctoas} to define the genuine function is given by the equation
\[
(\theta_{\pm}^2 \circ N_{F_{\alpha}/F} \circ \alpha^{\vee})(1+z)=\Lambda (\tr_{F_{\alpha}/F}(2\overline{b}_{\alpha}z))
\]
for \(\alpha \in R_{\sym}\) such that the depth of \(\theta_{\pm}^2 \circ N_{F_{\alpha}/F} \circ \alpha^{\vee}\) is positive.
Here we write \(\overline{b}_{\alpha}\) to distinguish it from \(\overline{a}_{\alpha}\).
For \(\alpha \in R_{\sym}\), we have the following commutative diagram by Proposition \ref{prop:mapJOS} \ref{prop:mapJOSsym}:
\[
\begin{tikzcd}
& F_{\alpha}^{\times}/N_{F_{\alpha}/F_{\pm \alpha}}(F_{\alpha}^{\times}) & S(F)_{\pm} & \bbC^{\times} \\
F_{\alpha}^{\times} & F_{\alpha}^1 & S(F) &
\arrow["\widetilde{N}_{\alpha}", twoheadrightarrow, from=2-1, to=2-2]
\arrow[twoheadrightarrow, bend left=20, from=2-1, to=1-2]
\arrow["\widetilde{N}_{\alpha}"', twoheadrightarrow, from=1-2, to=2-2]
\arrow["\mapJOS{\alpha, \pm}", from=1-2, to=1-3]
\arrow[twoheadrightarrow, from=1-3, to=2-3]
\arrow["\theta_{\pm}^2", from=1-3, to=1-4]
\arrow["N_{F_{\alpha}/F} \circ \alpha^{\vee}", bend right=20, from=2-1, to=2-3]
\arrow["\mapJOS{\alpha}", from=2-2, to=2-3]
\arrow[dashed, "\theta_{\pm}^2"', from=2-3, to=1-4]
\end{tikzcd}
\]
This shows that
\(
\theta_{\pm}^2 \circ N_{F_{\alpha}/F} \circ \alpha^{\vee}=(\theta_{\pm, \alpha})^2.
\)
As \(p\) is odd, the depth of \((\theta_{\pm, \alpha})^2\) is the same as that of \(\theta_{\pm, \alpha}\).
Hence, the depth of \(\theta_{\pm}^2\circ N_{F_{\alpha}/F} \circ \alpha^{\vee}\) is positive if and only if \(l_{\alpha}\) is positive,
and then we have \(2\overline{a}_{\alpha}=2\overline{b}_{\alpha}\) in \(\mathfrak{p}_{F_{\alpha}}^{-l_{\alpha}+e_{\alpha}r_{\alpha}} / \mathfrak{p}_{F_{\alpha}}^{-l_{\alpha}+e_{\alpha}r_{\alpha}+1}\).
Since \(2\neq p\) is invertible in \(\mathfrak{p}_{F_{\alpha}}^{-l_{\alpha}+e_{\alpha}r_{\alpha}} / \mathfrak{p}_{F_{\alpha}}^{-l_{\alpha}+e_{\alpha}r_{\alpha}+1}\), we deduce that \(\overline{a}_{\alpha}=\overline{b}_{\alpha}\).
For \(\alpha\) with \(l_{\alpha}=0\),
both constructions are independent of the choices of \(a_{\alpha}\), since \(F_{\alpha}/F_{\pm \alpha}\) is unramified.
Now we conclude that the resulting genuine functions are the same.
\end{proof}


\subsubsection{Relation with the canonical splitting of the cover} \label{sec:relspl}
We continue to assume that \(\theta_{\pm}\) satisfies {\bf (Ram)}.
For \(\alpha \in R_{\sym}\), another natural way to obtain a character of (an open subgroup of) \(F_{\alpha}^{\times}\) from \(\theta_{\pm}\) is to
descend a suitable restriction of \(\theta_{\pm} \circ \mapJOS{\alpha, \pm}\) to an open subgroup of \(J_{O}(F)=F_{\alpha}^1\) via the canonical splitting (cf.\ Proposition \ref{prop:minramchi}) and compose \(\widetilde{N}_{\alpha} \colon F_{\alpha}^{\times} \twoheadrightarrow F_{\alpha}^1=J_{O}(F)\).
Let us prove that the resulting character is essentially the restriction of \(\theta_{\pm, \alpha}\) and thus gives rise to the same genuine function.

To describe the canonical splitting, we take a set of \(\chi\)-data \(\chi=\{\chi_{\alpha}\}_{\alpha \in R(S, G)}\).
We assume
\begin{itemize}
	\item \(\chi_{\alpha}=1\) if \(\alpha \in R_{\asym}\);
	\item \(\chi_{\alpha}\) is unramified if \(\alpha \in R_{\sym}\) and \(F_{\alpha}/F_{\pm \alpha}\) is unramified;
	\item \(\chi_{\alpha}\) is tamely ramified if \(\alpha \in R_{\sym}\) and \(F_{\alpha}/F_{\pm \alpha}\) is ramified with \(d_{\alpha}=0\);
	\item \(\chi_{\alpha}\) is of level \(2d_{\alpha}\) as in Proposition \ref{prop:minramchi} if \(\alpha \in R_{\sym}\) and \(F_{\alpha}/F_{\pm \alpha}\) is ramified with \(d_{\alpha}>0\).
\end{itemize}
Note that if \(d_{\alpha}=0\) for every \(\alpha \in R_{\sym}\), then this is nothing but assuming that \(\chi\) is minimally ramified in the sense of \cite[Definition 4.6.1]{MR4013740}.

For \(\alpha \in R_{\sym}\) such that \(\theta_{\pm} \circ \mapJOS{\alpha, \pm} \colon F_{\alpha}^{\times}/N_{F_{\alpha}/F_{\pm \alpha}}(F_{\alpha}^{\times})\to \bbC ^{\times}\) is a genuine character, 
the product \(\left( \theta_{\pm} \circ \mapJOS{\alpha, \pm}\right)\cdot \chi_{\alpha}^{-1}\) descends to a character of \(F_{\alpha}^1\).
For \(\alpha \in R_{\sym}\) such that \(\theta_{\pm} \circ \mapJOS{\alpha, \pm}\) is not genuine, \(\theta_{\pm} \circ \mapJOS{\alpha, \pm}\) descends to a character of \(F_{\alpha}^1\).
In both cases, we then compose \(\widetilde{N}_{\alpha} \colon F_{\alpha}^{\times} \twoheadrightarrow F_{\alpha}^1\) to get a character \(\theta_{\pm, \alpha}' \colon F_{\alpha}^{\times}\to \bbC ^{\times}\).
The following commutative diagram summarizes the situation:
\[
\begin{tikzcd}
	& F_{\alpha}^{\times}/N_{F_{\alpha}/F_{\pm \alpha}}(F_{\alpha}^{\times}) & \bbC^{\times} \\
	F_{\alpha}^{\times} & F_{\alpha}^1 &
	\arrow[dashed, from=2-2, to=1-3]
	\arrow["\widetilde{N}_{\alpha}", twoheadrightarrow, from=2-1, to=2-2]
	\arrow[twoheadrightarrow, bend left=20, from=2-1, to=1-2]
	\arrow["\widetilde{N}_{\alpha}"', twoheadrightarrow, from=1-2, to=2-2]
	\arrow["(\ast)", from=1-2, to=1-3]
	\arrow["\theta_{\pm, \alpha}'"', bend right=45, from=2-1, to=1-3]
\end{tikzcd}
\]
Here, \((\ast)\) denotes 
\(\left( \theta_{\pm} \circ \mapJOS{\alpha, \pm}\right)\cdot \chi_{\alpha}^{-1}\) or \(\theta_{\pm} \circ \mapJOS{\alpha, \pm}\), and the homomorphism \(F_{\alpha}^{\times} \twoheadrightarrow F_{\alpha}^{\times}/N_{F_{\alpha}/F_{\pm \alpha}}(F_{\alpha}^{\times})\) is the canonical surjection.
Note that \(\theta_{\pm, \alpha}'\) can also be regarded as the composite of \((\ast)\) and the canonical surjection, and hence the restriction of \(\theta_{\pm, \alpha}'\) to \(U_{F_{\alpha}}^{2d_{\alpha}}\) (resp.\ to \(U_{F_{\alpha}}\)) is independent of the choice of \(\chi\) as above if \(d_{\alpha}>0\) by Proposition \ref{prop:minramchi} (resp.\ by the discussion following Proposition \ref{prop:minramchi}).

Recall that \(l_{\alpha}\) is the level of \(\theta_{\pm, \alpha}\) and put \(i_{\alpha}=l_{\alpha}-2d_{\alpha}\) as in the proof of Theorem \ref{thm:main1}.
If \(l_{\alpha}>0\), then 
\[
\theta_{\pm, \alpha}|_{U_{F_{\alpha}}^{i_{\alpha}}}=\theta_{\pm, \alpha}'|_{U_{F_{\alpha}}^{i_{\alpha}}}
\]
because \(\chi_{\alpha}|_{U_{F_{\alpha}}^{i_{\alpha}}}\) is trivial; for \(\alpha\) with \(d_{\alpha}>0\), this can be checked as
\[
i_{\alpha}-(2d_{\alpha}+1)
=(l_{\alpha}-2d_{\alpha})-(2d_{\alpha}+1)
=l_{\alpha}-4d_{\alpha}-1\geq 0
\]
by {\bf (Ram\(_{\alpha}\))}.
Therefore, even if we replace \(\theta_{\pm, \alpha}\) by \(\theta_{\pm, \alpha}'\) in \eqref{eq:linearization4},
the resulting genuine function will be the same.

\begin{rem} \label{rem:epsconstr}
For a symmetric root \(\alpha\) with \(l_{\alpha}>0\), the factor \(\kappa_{\alpha}\left(\frac{\delta_{\alpha}-\delta_{-\alpha}}{a_{\alpha}}\right)\) in \eqref{eq:aSfroma} admits an alternative description in terms of the character \(\chi_{\alpha}\) as above and local epsilon factors as follows.

Fix such an \(\alpha \in R_{\sym}\)
and put \(t_{\alpha}:=\left\lfloor \frac{l_{\alpha}}{2}\right\rfloor +1\). Then, as in Theorem \ref{thm:main1} \ref{item:linearization1},
there exists a unique element
\[
\overline{c}_{\alpha} \in \frakp_{F_{\alpha}}^{-l_{\alpha}+e_{\alpha}r_{\alpha}-d_{\alpha}}/\frakp_{F_{\alpha}}^{-t_{\alpha}+e_{\alpha}r_{\alpha}-d_{\alpha}+1}
\]
such that
\[
\theta_{\pm, \alpha}(1+z)=\Lambda (\tr_{F_{\alpha}/F}(\overline{c}_{\alpha}z))
\]
holds for any \(z\in \frakp_{F_{\alpha}}^{t_{\alpha}}\).
Comparing the defining conditions with \(\overline{a}_{\alpha}\) in Theorem \ref{thm:main1} \ref{item:linearization1}, we see that \(a_{\alpha}-c_{\alpha} \in \frakp_{F_{\alpha}}^{-l_{\alpha}+e_{\alpha}r_{\alpha}+d_{\alpha}+1}\),
so that \(a_{\alpha}c_{\alpha}^{-1}\in U_{F_{\alpha}}^{2d_{\alpha}+1}\), and therefore 
\[
\chi_{\alpha}(c_{\alpha})=\chi_{\alpha}(a_{\alpha})
\]
for any lifts \(a_{\alpha}, c_{\alpha} \in \frakp_{F_{\alpha}}^{-l_{\alpha}+e_{\alpha}r_{\alpha}-d_{\alpha}}\) of \(\overline{a}_{\alpha}, \overline{c}_{\alpha}\).
Now, the conductor/level condition in Deligne's twist formula \cite{MR349635}*{Lemme 4.16} (see also \cite{MR2234120}*{23.8 Stability Theorem} for the statement in terms of levels)
is satisfied by
\(\chi_{\alpha}\) and \(\theta_{\pm \alpha}\).
Indeed, in the ramified case, \(\chi_{\alpha}\) has level \(2d_{\alpha}\) and {\bf (Ram\(_\alpha\))} gives \(l_{\alpha}>4d_{\alpha}\);
in the unramified case, \(\chi_{\alpha}\) is unramified and \(l_{\alpha}>0\).
Hence, the twist formula gives
\[
\frac{\epsilon (\theta_{\pm \alpha}\chi_{\alpha}, s, \Lambda_{F_{\alpha}})}{\epsilon (\theta_{\pm \alpha}, s, \Lambda_{F_{\alpha}})}
=\chi_{\alpha}^{-1}(c_{\alpha})
=\chi_{\alpha}^{-1}(a_{\alpha}),
\]
where \(\Lambda_{F_{\alpha}}:=\Lambda \circ \tr_{F_{\alpha}/F}\).

Consider the factor
\(\kappa_{\alpha}\left(\frac{\delta_{\alpha}-\delta_{-\alpha}}{a_{\alpha}}\right)\)
in \eqref{eq:aSfroma} corresponding to \(\alpha (\gamma)\neq 1\).
Put \(D_{\alpha}:=\delta_{\alpha}-\delta_{-\alpha}\neq 0\).
Since \(\frac{D_{\alpha}}{a_{\alpha}} \in F_{\pm \alpha}^{\times}\) and \(\chi_{\alpha}|_{F_{\pm \alpha}^{\times}}=\kappa_{\alpha}\), it follows that
\[
\kappa_{\alpha}\left(\frac{D_{\alpha}}{a_{\alpha}}\right)
=\chi_{\alpha}(D_{\alpha}) \cdot \chi_{\alpha}(a_{\alpha})^{-1}
=\chi_{\alpha}(D_{\alpha}) \cdot \frac{\epsilon (\theta_{\pm \alpha}\chi_{\alpha}, s, \Lambda_{F_{\alpha}})}{\epsilon (\theta_{\pm \alpha}, s, \Lambda_{F_{\alpha}})}.
\]
Equivalently, putting \(\Lambda_{F_{\alpha}, D_{\alpha}}(x):=\Lambda_{F_{\alpha}}(D_{\alpha}\cdot x)\), we can also write
\[
\kappa_{\alpha}\left(\frac{D_{\alpha}}{a_{\alpha}}\right)
=\frac{\epsilon (\theta_{\pm \alpha}\chi_{\alpha}, s, \Lambda_{F_{\alpha}, D_{\alpha}})}{\epsilon (\theta_{\pm \alpha}, s, \Lambda_{F_{\alpha}, D_{\alpha}})}.
\]

Note that, when \(F_{\alpha}/F_{\pm \alpha}\) is ramified, the inequality in {\bf (Ram\(_\alpha\))} is precisely the conductor/level condition required for the twist formula.
\end{rem}

\subsection{Method 2: Another construction using the exponential map}\label{sec:method2}
We continue with the double cover \(S(F)_{\pm}\), a genuine character \(\theta_{\pm}\) and a non-trivial additive character \(\Lambda\) of \(F\). 

Assume now that the characteristic of \(F\) is zero.
We recall from \cite{1907.05173}*{\S3.7} the exponential map for the double cover in the form needed below, and use it to give another construction of the genuine function \(a_S\).

For any finite extension \(L/F\), define a lattice \(\omega_L \subset L\) by
\[
\omega_L:=\left\{x\in L \; \middle| \; v_L(x) >\frac{v_L(p)}{p-1}\right\}.
\]
Then the usual exponential power series \(\exp(x)=\sum_{n=0}^{\infty} \frac{x^n}{n!}\) converges for \(x\in \omega_L \subset L\).
Let \(E/F\) be a finite Galois extension splitting \(S\).
Then the Lie algebra \(\Lie(S)(F)\) of \(S\) is naturally identified with \((X_{\ast}(S)\otimes E)^{\Gamma_{E/F}}\).
The exponential map for \(S(E)\) is expressed as
\(\exp:=\id \otimes \exp \colon \Lie(S)(\omega_E):=X_{\ast}(S)\otimes \omega_E \to S(E)\)
and restricts to the exponential map
\[
\exp \colon \Lie(S)(\omega_F):=(\Lie(S)(\omega_E))^{\Gamma_{E/F}} \to S(F)
\]
for \(S(F)\) by the \(\Gamma_{E/F}\)-equivariance.
For any finite extension \(L/F\), we use \(\Lie(S)(\omega_L)\) for the analogous domain of the exponential map for \(S(L)\).
Next, we put
\[
\Omega:=\left\{X\in \Lie(S)(\omega_F)\; \middle| \; \ v_{F_\alpha}(d\alpha(X)) > v_{F_\alpha}(2)+\frac{v_{F_\alpha}(p)}{p-1} \text{ for any \(\alpha \in R_{\sym}\)}\right\}.
\]
If \(p\neq 2\), then \(v_{F_{\alpha}}(2)=0\), so the additional condition is automatic and \(\Omega=\Lie(S)(\omega_F)\).
For \(p=2\), the condition ensures the convergence of \(\exp\left(d{\alpha}\left(X\right)/2\right)\).
The exponential map for \(S(F)_{\pm \pm}\) is expressed as 
\[
\exp_{\pm \pm} \colon \Omega \to S(F)_{\pm \pm}, \quad 
X \mapsto \left(\exp X,\left(\exp\left(d{\alpha}\left(X\right)/2\right)\right)_{\alpha}, (1_O)_O\right),
\]
where the element of \(S(F)_{\pm \pm}\) is expressed as in \eqref{eq:eltofdc} and \(1_O\) denotes the unit in \(\{\pm 1\} \subset J_O(F)_{\pm}\) for an asymmetric \(\Sigma\)-orbit \(O\subset R\).
We write \(\exp_{\pm}\) for the composite of \(\exp_{\pm \pm}\) and the surjection \(S(F)_{\pm \pm} \twoheadrightarrow S(F)_{\pm}\).

Let us construct the genuine function using the exponential map.
The relevant objects are summarized in the following diagram:
\begin{equation} \label{eq:expdiagram}
\begin{tikzcd}
& & S(F)_{\pm} \arrow[r, "\theta_{\pm}"] \arrow[d, twoheadrightarrow, "\pr_1" '] & \bbC^{\times} \\
F_{\alpha}^{\times} \arrow[r, "\alpha^{\vee}"] & S(F_{\alpha}) \arrow[r, "N_{F_{\alpha}/F}"] & S(F) & \\
F_{\alpha} \arrow[u, dotted, "\exp"] \arrow[r, "d\alpha^{\vee}"] & \Lie(S)(F_{\alpha}) \arrow[u, dotted, "\exp"] \arrow[r, "\tr_{F_{\alpha}/F}"] & \Lie(S)(F) \arrow[u, dotted, "\exp"] \arrow[uu, bend right=60, dotted, "\exp_{\pm}" '], &
\end{tikzcd}
\end{equation}
where \(\alpha \in R\) is a root as usual, 
\(\tr_{F_{\alpha}/F}\) denotes the differential of the norm map \(N_{F_{\alpha}/F}\),
and a dotted arrow means that the homomorphism is defined on a suitable open subgroup of the domain.
For the commutativity of the squares and the ``triangle", it is easy to show the following:
\begin{itemize}
\item if \(z\in F_{\alpha}\) lies in \(\omega_{F_{\alpha}}\), then \(\exp(z) \in F_{\alpha}^{\times}\) and \(\exp(d\alpha^{\vee}(z)) \in S(F_{\alpha})\) are both defined and we have \(\alpha^{\vee}(\exp(z))=\exp(d\alpha^{\vee}(z))\);
\item if \(X\in \Lie(S)(F_{\alpha})\) lies in \(\Lie(S)(\omega_{F_{\alpha}})\), then \(\exp(X) \in S(F_{\alpha})\) and \(\exp(\tr_{F_{\alpha}/F}(X)) \in S(F)\) are both defined and we have \(N_{F_{\alpha}/F}(\exp(X))=\exp(\tr_{F_{\alpha}/F}(X))\);
\item if \(Y\in \Lie(S)(F)\) lies in \(\Omega\), then \(\exp(Y) \in S(F)\) and \(\exp_{\pm}(Y) \in S(F)_{\pm}\) are both defined and we have \(\pr_1(\exp_{\pm}(Y)) = \exp(Y)\).
\end{itemize}
We retain the notation \(e_{\alpha}, d_{\alpha}, r_{\alpha}\) from \S \ref{sec:method1}.
In addition, for \(\alpha \in R\), we put
\begin{align*}
m_{\alpha}&:=\min \left\{ m\geq 0 \; \middle| \; m>\frac{v_{F_{\alpha}}(p)}{p-1} \text{ and } \tr_{F_{\alpha}/F}(d\alpha^{\vee}(\frakp_{F_{\alpha}}^m)) \subset \Omega \right\}, \\
l_{\alpha}^{\mathrm{exp}}&:=\min \left\{ l\geq m_{\alpha}-1 \; \middle| \; \theta_{\pm} \circ \exp_{\pm} \circ \tr_{F_{\alpha}/F} \circ d\alpha^{\vee} \text{ is trivial on } \frakp_{F_{\alpha}}^{l+1} \right\}.
\end{align*}
Note that, on \(\frakp_{F_{\alpha}}^{m_{\alpha}}\), the composite \(\exp_{\pm} \circ \tr_{F_{\alpha}/F} \circ d{\alpha}^{\vee}\) is defined and we have
\[
\pr_1 \circ \exp_{\pm} \circ \tr_{F_{\alpha}/F} \circ d{\alpha}^{\vee}
=N_{F_{\alpha}/F} \circ \alpha^{\vee} \circ \exp.
\]
Moreover, if \(p\neq 2\), then \(m_{\alpha}=\left[ \frac{v_{F_{\alpha}}(p)}{p-1} \right]+1\).

For \(\alpha \in R_{\sym}\), we consider the following condition:
\begin{description}
\item[(Exp\(_{\alpha}\))] 
\(l_{\alpha}^{\mathrm{exp}}-m_{\alpha} \geq 2d_{\alpha}\).
\end{description}
Note that, under the assumption {\bf (Exp\(_{\alpha}\))}, we never have \(l_{\alpha}^{\mathrm{exp}}=0\). 
We write {\bf (Exp)} for the condition that {\bf (Exp\(_{\alpha}\))} holds for every \(\alpha \in R_{\sym}\).

\begin{thm}
Let \(\theta_{\pm}\) be a genuine character and assume that it satisfies {\bf (Exp)}.

Then the following hold.
\begin{enumerate}
\item \label{item:linearization2}
For \(\alpha \in R_{\sym}\), 
the composite \(\exp_{\pm} \circ \tr_{F_{\alpha}/F} \circ d\alpha^{\vee}\) is defined on 
\(\frakp_{F_{\alpha}}^{l_{\alpha}^{\mathrm{exp}}-2d_{\alpha}}\), and 
there exists a unique element \(\overline{a}_{\alpha} \in \frakp_{F_{\alpha}}^{-l_{\alpha}^{\mathrm{exp}}+e_{\alpha}r_{\alpha}-d_{\alpha}} / \frakp_{F_{\alpha}}^{-l_{\alpha}^{\mathrm{exp}}+e_{\alpha}r_{\alpha}+d_{\alpha}+1}\) such that
\begin{equation} \label{eq:linearization5}
	\theta_{\pm} \circ \exp_{\pm} \circ \tr_{F_{\alpha}/F} \circ d\alpha^{\vee}(z) = \Lambda (\tr_{F_{\alpha}/F}(\overline{a}_{\alpha}z))
\end{equation}
for any \(z\in \frakp_{F_{\alpha}}^{l_{\alpha}^{\mathrm{exp}}-2d_{\alpha}}\).
\item \label{item:liftofa2}
For \(\alpha \in R_{\sym}\), the element \(\overline{a}_{\alpha}\) defined by \eqref{eq:linearization5} can be lifted to an element \(a_{\alpha} \in F_{\alpha}^0 \cap \frakp_{F_{\alpha}}^{-l_{\alpha}^{\mathrm{exp}}+e_{\alpha}r_{\alpha}-d_{\alpha}}\),
so that we can define a genuine function \(a_S\) by the formula \eqref{eq:aSfroma}.
\item \label{item:independence2}
The genuine function \(a_S\) is independent of the choice of the lift \(a_{\alpha} \in F_{\alpha}^0 \cap \frakp_{F_{\alpha}}^{-l_{\alpha}^{\mathrm{exp}}+e_{\alpha}r_{\alpha}-d_{\alpha}}\) for each \(\alpha \in R_{\sym}\).
\item \label{item:coincidence2}
If \(p\neq 2\), then the genuine function defined this way coincides with the one defined by the mod-\(a\)-data of \cite{1907.05173}*{\S4.3} recalled in \S \ref{sec:gctoas}.
\end{enumerate}
\end{thm}
\begin{proof}
For \(\alpha \in R_{\sym}\), we put
\[
i_{\alpha}:=l_{\alpha}^{\mathrm{exp}}-2d_{\alpha}, \quad j_{\alpha}:=l_{\alpha}^{\mathrm{exp}},
\]
so that \(-j_{\alpha}+e_{\alpha}r_{\alpha}-d_{\alpha}=-l_{\alpha}^{\mathrm{exp}}+e_{\alpha}r_{\alpha}-d_{\alpha}\) and
\(-i_{\alpha}+e_{\alpha}r_{\alpha}-d_{\alpha}+1=-l_{\alpha}^{\mathrm{exp}}+e_{\alpha}r_{\alpha}+d_{\alpha}+1\).

To prove the assertions \ref{item:linearization2}, \ref{item:liftofa2}, take \(\alpha \in R_{\sym}\).
We have \(i_{\alpha}\geq m_{\alpha}\) by {\bf (Exp\(_{\alpha}\))}, so that 
\(\theta_{\pm} \circ \exp_{\pm} \circ \tr_{F_{\alpha}/F} \circ d\alpha^{\vee}\) is defined on \(\frakp_{F_{\alpha}}^{i_{\alpha}}\) and the left-hand side of \eqref{eq:linearization5} makes sense.
Then Proposition \ref{prop:liftofa} \ref{prop:liftofa1} implies the assertion \ref{item:linearization2}.

If we write
\[
d\widetilde{N}_{\alpha} \colon F_{\alpha} \twoheadrightarrow F_{\alpha}^0, \quad z \mapsto z-\tau_{\alpha}(z)
\]
for the differential of \(\widetilde{N}_{\alpha} \colon F_{\alpha}^{\times} \twoheadrightarrow F_{\alpha}^1\), then we have \(\tr_{F_{\alpha}/F_{\pm \alpha}} \circ d\alpha^{\vee}=d\alpha^{\vee} \circ d\widetilde{N}_{\alpha} \colon F_{\alpha} \to \Lie(S)(F_{\pm \alpha})\).
This shows that
\(\theta_{\pm} \circ \exp_{\pm} \circ \tr_{F_{\alpha}/F} \circ d\alpha^{\vee}\) is trivial on \(\frakp_{F_{\alpha}}^{m_{\alpha}} \cap \Ker(d\widetilde{N}_{\alpha})=\frakp_{F_{\alpha}}^{m_{\alpha}} \cap F_{\pm \alpha}\).
It follows from Proposition \ref{prop:liftofa} \ref{prop:liftofa2} that \(\overline{a}_{\alpha}\) lifts to an element \(a_{\alpha}\in F_{\alpha}^0 \cap \frakp_{F_{\alpha}}^{-j_{\alpha}+e_{\alpha}r_{\alpha}-d_{\alpha}}\).

Let us prove the assertion \ref{item:independence2}.
For each \(\alpha \in R_{\sym}\),
since we have \(v_{F_{\alpha}}(a_{\alpha})=-j_{\alpha}+e_{\alpha}r_{\alpha}-d_{\alpha}\)
by the definition of \(l_{\alpha}^{\mathrm{exp}}\)
and
\[
(-i_{\alpha}+e_{\alpha}r_{\alpha}-d_{\alpha}+1)-(-j_{\alpha}+e_{\alpha}r_{\alpha}-d_{\alpha})
=j_{\alpha}-i_{\alpha}+1=2d_{\alpha}+1,
\]
Proposition \ref{prop:wildmoda} shows that \(\kappa_{\alpha}\left(\frac{\delta_{\alpha}-\delta_{-\alpha}}{a_{\alpha}}\right)\) in the formula \eqref{eq:aSfroma} is independent of the choice of the lift \(a_{\alpha}\).
Thus, \(a_S\) is independent of the choice of such lifts.

Let us finally prove the assertion \ref{item:coincidence2}.
Let \((\overline{b}_{\alpha})_{\alpha}\) be the mod-\(a\)-data used in \S \ref{sec:gctoas} to define the genuine function.
Let \(\alpha\) be a symmetric root.
Note that \(d_{\alpha}=0\) and thus \(i_{\alpha}=j_{\alpha}=l_{\alpha}^{\mathrm{exp}}\).
In view of the diagram
\[
\begin{tikzcd}
& & S(F)_{\pm} \arrow[r, "\theta_{\pm}^2"] \arrow[d, twoheadrightarrow, "\pr_1" '] & \bbC^{\times} \\
F_{\alpha}^{\times} \arrow[r, "\alpha^{\vee}"] & S(F_{\alpha}) \arrow[r, "N_{F_{\alpha}/F}"] & S(F) \arrow[ur, dashed, pos=0.6, "\theta_{\pm}^2"] & \\
F_{\alpha} \arrow[u, dotted, "\exp"] \arrow[r, "d\alpha^{\vee}"] & \Lie(S)(F_{\alpha}) \arrow[u, dotted, "\exp"] \arrow[r, "\tr_{F_{\alpha}/F}"] & \Lie(S)(F) \arrow[u, dotted, "\exp"] \arrow[uu, bend right=60, dotted, "\exp_{\pm}" '], &
\end{tikzcd}
\]
we have
\[
\theta_{\pm}^2 \circ \exp_{\pm} \circ \tr_{F_{\alpha}/F} \circ d{\alpha}^{\vee}
=\theta_{\pm}^2 \circ N_{F_{\alpha}/F} \circ \alpha^{\vee} \circ \exp
\]
on \(\frakp_{F_{\alpha}}^{l_{\alpha}^{\mathrm{exp}}}\) since \(l_{\alpha}^{\mathrm{exp}} \geq m_{\alpha}\) by {\bf (Exp\(_\alpha\))},
and they both induce a non-trivial character of \(\frakp_{F_{\alpha}}^{l_{\alpha}^{\mathrm{exp}}}/\frakp_{F_{\alpha}}^{l_{\alpha}^{\mathrm{exp}}+1}\) since \(p\) is odd.
It follows that the depth of \(\theta_{\pm}^2 \circ N_{F_{\alpha}/F} \circ \alpha^{\vee}\) is positive.
Hence, we have, for any \(z\in \frakp_{F_{\alpha}}^{l_{\alpha}^{\mathrm{exp}}}\),
\begin{align*}
\Lambda (\tr_{F_{\alpha}/F}(2\overline{a}_{\alpha}z))
&= \theta_{\pm}^2 \circ \exp_{\pm} \circ \tr_{F_{\alpha}/F} \circ d\alpha^{\vee}(z) \\
&=\theta_{\pm}^2 \circ N_{F_{\alpha}/F} \circ \alpha^{\vee}(\exp(z)) \\
&=\theta_{\pm}^2 \circ N_{F_{\alpha}/F} \circ \alpha^{\vee}(1+z) \\
&=\Lambda (\tr_{F_{\alpha}/F}(2\overline{b}_{\alpha}z)).
\end{align*} 
Here, for the third equality, put \(l=l_{\alpha}^{\mathrm{exp}}\). Since \(l\geq m_{\alpha}>v_{F_{\alpha}}(p)/(p-1)\), for every \(n\geq2\) we have

\[
v_{F_{\alpha}}\left(\frac{z^n}{n!}\right)
\geq nl-\frac{v_{F_{\alpha}}(p)(n-1)}{p-1}>l,
\]
so that
\[
\frac{\exp(z)}{1+z}\in U_{F_{\alpha}}^{l+1}
\subset
\Ker\left(\theta_{\pm}^2\circ N_{F_{\alpha}/F}\circ\alpha^\vee\right).
\]
It follows that \(2\overline{a}_{\alpha}=2\overline{b}_{\alpha}\) in \(\frakp_{F_{\alpha}}^{-l_{\alpha}^{\mathrm{exp}}+e_{\alpha}r_{\alpha}} / \frakp_{F_{\alpha}}^{-l_{\alpha}^{\mathrm{exp}}+e_{\alpha}r_{\alpha}+1}\).
Now we conclude that \(\overline{a}_{\alpha}=\overline{b}_{\alpha}\) by the invertibility of \(2\neq p\) in \(\frakp_{F_{\alpha}}^{-l_{\alpha}^{\mathrm{exp}}+e_{\alpha}r_{\alpha}} / \frakp_{F_{\alpha}}^{-l_{\alpha}^{\mathrm{exp}}+e_{\alpha}r_{\alpha}+1}\) and thus they give rise to the same genuine function.
\end{proof}

\section{Example: \(\PGL_2\) case} \label{sec:PGL2}
Here we illustrate our general theory in the case \(G=\PGL_2\).
Let \(S\) be an elliptic maximal torus.
Then \(S \simeq \Res_{E/F}^1 \bbG_m\) for a quadratic separable extension \(E/F\).
Let \(\theta_{\pm}\) be a genuine character of the double cover \(S(F)_{\pm}\) associated to the set \(R\subset X^{\ast}(S)\) of absolute roots of \(G\) with respect to \(S\).

Let us fix some symbols.
\begin{itemize}
\item We simply write \(\tau\) for the non-trivial element in \(\Gamma_{E/F}\) and \(\kappa\) for the isomorphism \(F^{\times}/N_{E/F}(E^{\times}) \stackrel{\sim}{\to} \{\pm 1\}\).
\item We identify \(S(E)=E^{\times}\), on which \(\tau\) acts as \(x\mapsto \tau(x)^{-1}\).
Thus, we have \(S(F)=E^1\), which is also identified with \(E^{\times}/F^{\times}\) by \(E^{\times}/F^{\times}\stackrel{\sim}{\to} E^1, \quad (x \mod F^{\times})\mapsto \frac{x}{\tau(x)}\).
Also, the norm map \(N_{E/F}\) for the torus \(S\) 
and its differential \(\tr_{E/F}\)
are given by 
\begin{align*}
&\widetilde{N}_{\alpha} \colon S(E)=E^{\times} \to S(F)=E^1, \quad x\mapsto \frac{x}{\tau(x)}, \\
&d\widetilde{N}_{\alpha} \colon \Lie(S)(E)=E \to \Lie(S)(F)=E^0, \quad x\mapsto x-\tau(x).
\end{align*}
\item Let \(\alpha \colon S \to \bbG_m\) be the $E$-homomorphism that induces the identity map on \(E\)-valued points: \(S(E)=E^{\times} \stackrel{\text{id}}{\to} E^{\times}=\bbG_m(E)\).
Then we have \(R=\{\alpha, -\alpha\}\) and \(\alpha\) is symmetric.
Moreover, the coroot \(\alpha^{\vee} \colon \bbG_m \to S\) is expressed on \(E\)-valued points as
\[
\alpha ^{\vee} \colon \bbG_m(E)=E^{\times} \to S(E)=E^{\times}, \quad x\mapsto x^2.
\]
\item With the choice of \(\alpha\), we can identify \(S(F)_{\pm}\) with \(E^{\times}/N_{E/F}(E^{\times})\) and express the covering map as
\[
\widetilde{N}_{\alpha} \colon S(F)_{\pm}=E^{\times}/N_{E/F}(E^{\times}) \to S(F)=E^1, \quad (x\mod N_{E/F}(E^{\times}))\mapsto \frac{x}{\tau(x)}.
\]
\item We simply write \(d\) for \(d_{\alpha}\), so that the relative discriminant \(\frakd_{E/F}\) equals \(\frakp_F^{d+1}\) if \(E/F\) is ramified.
\item Let \(l_0\) be the minimum of the levels of the character twists \(\theta_{\pm} \cdot (\eta \circ N_{E/F})\), where \(\theta_{\pm}\) is viewed as a character of \(E^{\times}\) and \(\eta\) ranges over all characters of \(F^{\times}\).
\end{itemize}
We assume that,
when \(\theta_{\pm}\) is identified with a character of the Weil group \(W_E\) of \(E\), the induction \(\Ind_{W_E}^{W_F}\theta_{\pm}\) is irreducible.
Combining the local Langlands correspondence for \(S(F)_{\pm}\) \cite{1907.05173}*{\S 3.4} 
with the canonical \(L\)-embedding \({}^LS_{\pm}\to{}^LG\) \cite{1907.05173}*{\S 4.1},
we obtain from \(\theta_{\pm}\) an \(L\)-parameter for \(G\).
This is the forward map appearing in \cite{1907.05173}*{Proposition 4.2.1}; unlike the bijectivity statement there, its construction does not require the assumption \(p\nmid |W|\).
In the present \(G=\PGL_2\) case, the explicit description in \cite{1907.05173}*{\S B.2},
together with a direct inspection of the canonical \(L\)-embedding, identifies this parameter with \(\Ind_{W_E}^{W_F}\theta_{\pm}\).
Note that, since \(\theta_{\pm}|_{F^{\times}}=\kappa\), this induced representation indeed has trivial determinant (see e.g.\ \cite[Proposition 29.2]{MR2234120})
and hence takes values in \(\widehat{G}=\SL_2(\mathbb C)\).
Thus, our irreducibility assumption says precisely that the resulting \(L\)-parameter is supercuspidal.
Moreover, by \cite[Theorem 41.5]{MR2234120}, a supercuspidal \(L\)-parameter for \(G\) normalizing a maximal torus arises in this way from a genuine character \(\theta_{\pm}\) satisfying \(l_0\geq 2d\).
Thus, we assume this inequality throughout this section.

\subsection{Method 1 for \(\PGL_2\)}
Recall that \(l_{\alpha}\) is the level of \(\theta_{\pm, \alpha}\).
\begin{prop} \label{prop:levelalpha}
We have
\[
l_{\alpha}=
l_0-d.
\]
\end{prop}
\begin{proof}
Note first that the definition \eqref{eq:defNbent} of \(\Nbent{\alpha}{\beta}\) for \(\beta=\alpha\) works out to
\begin{align*}
\Nbent{\alpha}{\alpha}\colon F_{\alpha}^{\times}=E^{\times} &\to F_{\alpha}^{\times}/N_{F_{\alpha}/F_{\pm \alpha}}(F_{\alpha}^{\times})=E^{\times}/N_{E/F}(E^{\times}) \\
x&\mapsto (x^2 \mod{N_{E/F}(E^{\times})})=\left(\frac{x}{\tau(x)} \mod{N_{E/F}(E^{\times})}\right).
\end{align*}
It follows that \(\theta_{\pm, \alpha}(x)=\theta_{\pm}\left(\frac{x}{\tau(x)}\right)=\theta_{\pm}(\widetilde{N}_{\alpha}(x))\).
Thus, if we express \(\theta_{\pm}=\theta_{\pm, 0} \cdot (\eta \circ N_{E/F})\) with a character \(\theta_{\pm, 0}\) of \(E^{\times}\) of level \(l_0\) and a character \(\eta\) of \(F^{\times}\), we have
\begin{equation} \label{eq:alphapart}
\theta_{\pm, \alpha}(x)
=\left(\theta_{\pm, 0} \cdot (\eta \circ N_{E/F})\right) \left(\widetilde{N}_{\alpha}(x)\right)
=\theta_{\pm, 0} \left(\widetilde{N}_{\alpha}(x)\right).
\end{equation}

Suppose that \(d>0\). Then \(l_0-d\) is odd by \cite[41.4 Lemma (2)]{MR2234120} and thus we have an isomorphism
\[
\widetilde{N}_{\alpha} \colon U_E^{l_0-d}/U_E^{l_0-d+1} \xrightarrow{\sim} U_E^{l_0}/U_E^{l_0+1}
\]
by \cite[41.2 Lemma]{MR2234120}.
Since the level of \(\theta_{\pm, 0}\) is \(l_0\), we see that \(l_{\alpha}=l_0-d\). 

It is similar but easier to see \(l_{\alpha}=l_0\) if \(d=0\).
\end{proof}
By Proposition \ref{prop:levelalpha}, the assumption 
{\bf (Ram\(_{\alpha}\))} in this case reads
\begin{description}
\item[(Ram\(_{\alpha}\))] if \(E/F\) is ramified, then \(l_0\geq 5d+1\).
\end{description}
\begin{rem}
Suppose that \(d>0\).
Let us briefly compare the inequality \(l_0\geq 5d+1\) above with an inequality appearing in the local Langlands correspondence for \(G\) explained in \cite{MR2234120}.
If \(l_0> 3d\), then, in the notation of \cite[Proposition 44.7]{MR2234120}, we have \(\delta_{E/F} \in \frakp^{-[m/2]}\) and thus \(F(\alpha)\simeq E=F(\beta)\).
This means that the unique (supercuspidal) representation \(\pi\) of \(G(F)\) in the \(L\)-packet of \(\Ind_{W_E}^{W_F}\theta_{\pm}\) is induced from a character of a compact open subgroup containing \(S(F)\),
rather than one containing \(S'(F)\) for a different elliptic maximal torus \(S'\subset G\).
The inequality \(l_0>3d\) is strictly weaker than \(l_0 \geq 5d+1\), and we do not know any natural interpretation of the latter stronger inequality in this context.
\end{rem}

\subsection{Method 2 for \(\PGL_2\)}
Assume that \(F\) is of characteristic zero
and let us explicate \(m_{\alpha}\), \(l_{\alpha}^{\mathrm{exp}}\) and {\bf (Exp\(_\alpha\))} in \S \ref{sec:method2} for \(G=\PGL_2\).
The diagram \eqref{eq:expdiagram} for \(\PGL_2\) reads as follows:
\[\begin{tikzcd}
& & E^{\times}/N_{E/F}(E^{\times}) \arrow[r, "\theta_{\pm}"] \arrow[d, twoheadrightarrow] & \bbC^{\times} \\
E^{\times} \arrow[r, "\alpha^{\vee}"] & E^{\times} \arrow[r, "\widetilde{N}_{\alpha}"] & E^1 & \\
E \arrow[u, dotted, "\exp"] \arrow[r, "d\alpha^{\vee}"] & E \arrow[u, dotted, "\exp"] \arrow[r, "d\widetilde{N}_{\alpha}"] & E^0 \arrow[u, dotted, "\exp"] \arrow[uu, bend right=60, dotted, "\exp_{\pm}" '], &
\end{tikzcd}\]
where \(d\alpha^{\vee}(z)=2z\) for \(z\in E\) and 
\(\exp_{\pm}(y)=\exp(y/2)\) for \(y\in E^0\).
In particular, in this case, if \(v_E(z) > \frac{v_E(p)}{p-1}\), then \(v_E(d\widetilde{N}_{\alpha}\circ d\alpha^{\vee}(z))=v_E(2(z-\tau(z))) > v_E(2) + \frac{v_E(p)}{p-1}\) holds automatically.
Thus, we have 
\begin{equation} \label{eq:malpha}
m_{\alpha}=\left[ \frac{v_E(p)}{p-1} \right]+1.
\end{equation}

\begin{prop} \label{prop:levelalpha2}
We have the following:
\begin{enumerate}
\item\label{item:levelalpha2} \(l_{\alpha}^{\mathrm{exp}}=\max \{m_{\alpha}-1, l_0-d\}\).
\item\label{item:expalpha} {\bf (Exp\(_\alpha\))} is equivalent to
\(l_0 \geq m_{\alpha}+3d\).
\end{enumerate}
\end{prop}
\begin{proof}
\ref{item:levelalpha2} 
To compute \(l_{\alpha}^{\mathrm{exp}}\), let us write \(\theta_{\pm}=\theta_{\pm, 0} \cdot (\eta \circ N_{E/F})\) as in the proof of Proposition \ref{prop:levelalpha}.
Then we have
\[
\theta_{\pm} \circ \exp_{\pm} \circ d\widetilde{N}_{\alpha} \circ d{\alpha}^{\vee}
=\theta_{\pm, 0} \circ \exp \circ d\widetilde{N}_{\alpha}
\]
on \(\frakp_E^{m_{\alpha}}\).
Indeed, for any \(z\in \frakp_E^{m_{\alpha}}\), we have
\[
\theta_{\pm} \circ \exp_{\pm} \circ d\widetilde{N}_{\alpha} \circ d{\alpha}^{\vee}(z)
=\theta_{\pm}(\exp(d\widetilde{N}_{\alpha}(z)))
=\theta_{\pm, 0}(\exp(d\widetilde{N}_{\alpha}(z)))
\]
as \(\exp(d\widetilde{N}_{\alpha}(z))=\widetilde{N}_{\alpha}(\exp(z)) \in E^1\).

We have \(d\widetilde{N}_{\alpha}(\frakp_{E}^{i})\subset \frakp_{E}^{i+d}\) for any \(i\in \bbZ\).
Therefore, if \(m_{\alpha}+d \geq l_0+1\), then \(\theta_{\pm} \circ \exp_{\pm} \circ d\widetilde{N}_{\alpha} \circ d\alpha^{\vee}\) is trivial on \(\frakp_E^{m_{\alpha}}\) and thus we have
\(
l_{\alpha}^{\mathrm{exp}}=m_{\alpha}-1
\)
in this case.

The computation of \(l_{\alpha}^{\mathrm{exp}}\) in the remaining case \(l_0+1 >m_{\alpha}+d\) is similar to the proof of Proposition \ref{prop:levelalpha}.
Suppose first that \(d>0\).
Then \(l_0-d\) is odd by \cite{MR2234120}*{41.4 Lemma (2)} and thus we have isomorphisms
\[
\frakp_{E}^{l_0-d}/\frakp_{E}^{l_0-d+1}\xrightarrow[d\widetilde{N}_{\alpha}]{\sim} \frakp_{E}^{l_0}/\frakp_{E}^{l_0+1}
\xrightarrow[\exp]{\sim} U_E^{l_0}/U_E^{l_0+1},
\]
where the first isomorphism follows from an additive analogue of \cite{MR2234120}*{41.2 Lemma}
and \(\exp\) is defined on \(\frakp_E^{l_0}\) as \(l_0\geq m_{\alpha}\) by \(d>0\).
Since \(\theta_{\pm, 0}\) has level \(l_0\), we see that \(l_{\alpha}^{\mathrm{exp}}=l_0-d\).
It is similar but easier to show \(l_{\alpha}^{\mathrm{exp}}=l_0\) if \(d=0\).

To summarize, we have
\[
l_{\alpha}^{\mathrm{exp}}=\begin{cases}
m_{\alpha}-1 &\text{ if \(m_{\alpha}+d\geq l_0+1\)} \\
l_0-d &\text{ otherwise}
\end{cases}
=\max \{m_{\alpha}-1, l_0-d\}.
\]

\ref{item:expalpha}
If \(m_{\alpha}+d\geq l_0+1\), then \(l_{\alpha}^{\mathrm{exp}}=m_{\alpha}-1\) by \ref{item:levelalpha2}, and the assumption {\bf (Exp\(_\alpha\))} is never satisfied.
Thus, the assumption {\bf (Exp\(_\alpha\))} amounts to \(l_0+1>m_{\alpha}+d\) and \((l_0-d)-m_{\alpha} \geq 2d\).
As the former condition follows from the latter, this is equivalent to
\(l_0\geq m_{\alpha}+3d\), as stated.
\end{proof}
\begin{rem}
When \(d>0\), {\bf (Exp\(_\alpha\))} is neither stronger nor weaker than {\bf (Ram\(_\alpha\))} in general.
Indeed, by \eqref{eq:malpha} and Proposition \ref{prop:levelalpha2} \ref{item:expalpha}, the assumption {\bf (Exp\(_\alpha\))} reads
\[
l_0 \geq 3d+v_E(2)+1.
\]
We have \(1\leq d\leq v_E(2)\) (see for example \cite[41.1 Lemma (2)]{MR2234120}), but 
\[
(3d+v_E(2)+1)-(5d+1)=v_E(2)-2d
\]
can be either positive, zero, or negative.
\end{rem}

\bibliography{bibliography}

@misc{1912.03274,
Author = {Tasho Kaletha},
Title = {Supercuspidal {L}-packets},
Year = {2019},
Eprint = {arXiv:1912.03274},
}

@incollection{MR349635,
	author = {Deligne, P.},
	booktitle = {Modular functions of one variable, {II} ({P}roc. {I}nternat. {S}ummer {S}chool, {U}niv. {A}ntwerp, {A}ntwerp, 1972)},
	mrclass = {12A70 (10H10 12A65 12B25 12B30)},
	mrnumber = {349635},
	mrreviewer = {Roger\ Howe},
	pages = {501--597},
	publisher = {Springer, Berlin-New York},
	series = {Lecture Notes in Math.},
	title = {Les constantes des \'equations fonctionnelles des fonctions {$L$}},
	volume = {Vol. 349},
	year = {1973}}

@article{1907.05173,
	author = {Kaletha, Tasho},
	doi = {10.1353/ajm.2026.a1002870},
	fjournal = {American Journal of Mathematics},
	journal = {Amer. J. Math.},
	number = {5},
	pages = {1229--1300},
	title = {On {L}-embeddings and double covers of tori over local fields},
	url = {https://doi.org/10.1353/ajm.2026.a1002870},
	volume = {148},
	year = {2026}}

@article{MR909227,
	author = {Langlands, R. P. and Shelstad, D.},
	doi = {10.1007/BF01458070},
	fjournal = {Mathematische Annalen},
	issn = {0025-5831},
	journal = {Math. Ann.},
	mrclass = {11S37 (11F70 11F72 22E35 22E45 22E50)},
	mrnumber = {909227},
	mrreviewer = {Ernst-Wilhelm Zink},
	number = {1-4},
	pages = {219--271},
	title = {On the definition of transfer factors},
	url = {https://doi.org/10.1007/BF01458070},
	volume = {278},
	year = {1987}}

@book{MR2234120,
	author = {Bushnell, Colin J. and Henniart, Guy},
	doi = {10.1007/3-540-31511-X},
	isbn = {978-3-540-31486-8; 3-540-31486-5},
	mrclass = {22E50 (11-02 11S37 22-02)},
	mrnumber = {2234120},
	mrreviewer = {Alexandru Ioan Badulescu},
	pages = {xii+347},
	publisher = {Springer-Verlag, Berlin},
	series = {Grundlehren der Mathematischen Wissenschaften [Fundamental Principles of Mathematical Sciences]},
	title = {The local {L}anglands conjecture for {$\rm GL(2)$}},
	url = {https://doi.org/10.1007/3-540-31511-X},
	volume = {335},
	year = {2006}}

@article{MR4013740,
	author = {Kaletha, Tasho},
	doi = {10.1090/jams/925},
	fjournal = {Journal of the American Mathematical Society},
	issn = {0894-0347},
	journal = {J. Amer. Math. Soc.},
	mrclass = {22E50 (11F70 11S37)},
	mrnumber = {4013740},
	number = {4},
	pages = {1071--1170},
	title = {Regular supercuspidal representations},
	url = {https://doi.org/10.1090/jams/925},
	volume = {32},
	year = {2019}}

\noindent
L\'eo Gratien\\
Department of Mathematics, University of Chicago, Chicago, IL 60637\\
gratien@math.uchicago.edu \vspace*{10pt}

\noindent
Nhat Hoang Le\\
Department of Mathematics, Aarhus University, Ny Munkegade 118, 8000 Aarhus, Denmark\\
le.nhathoang@math.au.dk \vspace*{10pt}

\noindent
Kazuki Tokimoto\\
Faculty of Innovative Information Science, Chiba Institute of Technology, 2-1-1 Shibazono, Narashino, Chiba 275-0023, Japan\\
tokimoto.kazuki@chibatech.ac.jp

\end{document}